\documentclass[a4paper,11pt]{article}
\usepackage{amsmath,amsthm}
\usepackage{authblk}
\usepackage[style=numeric-comp,maxbibnames=99,bibencoding=utf8,firstinits=true]{biblatex}
\usepackage{booktabs}
\usepackage{cancel}
\usepackage{colortbl}
\usepackage[nomarkers,nolists]{endfloat}
\usepackage[hmargin={25mm,25mm},vmargin={30mm,35mm}]{geometry}
\usepackage{loglogslopetriangle}
\usepackage{multirow}
\usepackage{nicefrac}
\usepackage{paralist}
\usepackage[colorlinks, allcolors=blue]{hyperref}
\usepackage{pgfplots}
\usepackage{subcaption}
\usepackage[compact,small]{titlesec}
\usepackage{tikz}
\usepackage{tikz-cd}
\usepackage{newtxmath}
\usepackage{newtxtext}

\bibliography{curlcurl2d}
\newcommand{\email}[1]{\href{mailto:#1}{#1}}

\allowdisplaybreaks

\usepackage[normalem]{ulem}
\newcounter{corr}
\definecolor{violet}{rgb}{0.580,0.,0.827}
\newcommand{\corr}[3]{\typeout{Warning : a correction remains in page \thepage}
	\stepcounter{corr}        
				      {\color{blue}\ifmmode\text{\,\sout{\ensuremath{#1}}\,}\else\sout{#1}\fi}
             {\color{red}#2}
             {\color{violet} #3}
}

\newtheorem{theorem}{Theorem}

\newtheorem{lemma}[theorem]{Lemma}

\theoremstyle{remark}
\newtheorem{remark}[theorem]{Remark}
\theoremstyle{definition}

\newcommand{\st}{\,:\,}
\newcommand{\Real}{\mathbb{R}}

\DeclareRobustCommand{\bvec}[1]{\boldsymbol{#1}}
\newcommand{\uvec}[1]{\underline{\bvec{#1}}}

\newcommand{\cvec}[1]{\bvec{\mathcal{#1}}}

\DeclareRobustCommand{\btens}[1]{\boldsymbol{#1}}
\newcommand{\utens}[1]{\underline{\bvec{#1}}}

\DeclareMathOperator{\card}{card}

\newcommand{\rotation}[1]{\varrho_{#1}}

\DeclareMathOperator{\GRAD}{\bf grad}
\DeclareMathOperator{\DIV}{div}
\DeclareMathOperator{\ROT}{rot}
\DeclareMathOperator{\VROT}{\bf rot}

\newcommand{\Hrotrot}[1]{\bvec{H}(\VROT\ROT;#1)}
\newcommand{\Hrot}[1]{\bvec{H}(\ROT;#1)}

\newcommand{\HrotrotZ}[1]{\bvec{H}_0(\VROT\ROT;#1)}

\newcommand{\compl}{{\rm c}}

\newcommand{\Poly}[1]{\mathcal{P}^{#1}}
\newcommand{\vPoly}[1]{\cvec{P}^{#1}}
\newcommand{\Roly}[1]{\cvec{R}^{#1}}
\newcommand{\cRoly}[1]{\cvec{R}^{\compl,#1}}
\newcommand{\Goly}[1]{\cvec{G}^{#1}}
\newcommand{\cGoly}[1]{\cvec{G}^{\compl,#1}}

\newcommand{\lproj}[2]{\pi_{\mathcal{P},#2}^{#1}}
\newcommand{\vlproj}[2]{\bvec{\pi}_{\cvec{P},#2}^{#1}}
\newcommand{\Rproj}[2]{\bvec{\pi}_{\cvec{R},#2}^{#1}}
\newcommand{\cRproj}[2]{\bvec{\pi}_{\cvec{R},#2}^{\compl,#1}}

\newcommand{\edges}[1]{\mathcal{E}_{#1}}
\newcommand{\vertices}[1]{\mathcal{V}_{#1}}

\newcommand{\ET}{\edges{T}}

\newcommand{\VE}{\vertices{E}}
\newcommand{\VT}{\vertices{T}}

\newcommand{\normal}{\bvec{n}}
\newcommand{\tangent}{\bvec{t}}

\newcommand{\Mh}{\mathcal{M}_h}
\newcommand{\Th}{\mathcal{T}_h}
\newcommand{\Eh}{\mathcal{E}_h}
\newcommand{\Vh}{\mathcal{V}_h}

\newcommand{\boundary}{{\rm b}}
\newcommand{\Ehb}{\mathcal{E}_h^\boundary}
\newcommand{\Vhb}{\mathcal{V}_h^\boundary}

\DeclareMathOperator{\Ker}{Ker}
\DeclareMathOperator{\Image}{Im}

\newcommand{\norm}[2][]{\|#2\|_{#1}}
\newcommand{\seminorm}[2][]{|#2|_{#1}}
\newcommand{\vvvert}{\vert\kern-0.25ex\vert\kern-0.25ex\vert}
\newcommand{\tnorm}[2][]{\vvvert #2\vvvert_{#1}}

\newcommand{\term}{\mathfrak{T}}

\newcommand{\Xgrad}[2]{\underline{X}_{\GRAD,#2}^{#1}}
\newcommand{\Xrot}[1]{\uvec{X}_{\ROT,#1}^k}
\newcommand{\XrotZ}[1]{\uvec{X}_{\ROT,#1,0}^k}

\newcommand{\SV}[1]{\underline{V}_{#1}^k}
\newcommand{\SSigma}[1]{\uvec{\Sigma}_{#1}^k}
\newcommand{\SW}[1]{\underline{W}_{#1}^k}

\newcommand{\SVser}[1]{\widehat{\underline{V}}_{#1}^k}
\newcommand{\SSigmaser}[1]{\widehat{\uvec{\Sigma}}_{#1}^k}

\newcommand{\SVZ}[1]{\underline{V}_{#1,0}^k}
\newcommand{\SSigmaZ}[1]{\uvec{\Sigma}_{#1,0}^k}
\newcommand{\SWZ}[1]{\underline{W}_{#1,0}^k}

\newcommand{\Injrot}{\uvec{\mathfrak I}_{\ROT,h}^k}
\newcommand{\Resrot}{\uvec{\mathfrak S}_{\ROT,h}^k}
\newcommand{\ResrotT}{\uvec{\mathfrak S}_{\ROT,T}^k}

\newcommand{\IV}[1]{\underline{I}_{V,#1}^{k}}
\newcommand{\ISigma}[1]{\utens{I}_{\bvec{\Sigma},#1}^k}
\newcommand{\IW}[1]{\underline{I}_{W,#1}^k}
\newcommand{\Irot}[1]{\uvec{I}_{\ROT,#1}^k}

\newcommand{\GE}{G_E^k}
\newcommand{\GT}{\bvec{G}_T^k}
\newcommand{\PVT}{P_{V,T}^{k+1}}
\newcommand{\PVE}{P_{V,E}^{k+1}}
\newcommand{\PWT}{P_{W,T}^{k+1}}
\newcommand{\PWE}{P_{W,E}^{k+1}}
\newcommand{\PWF}{P_{W,F}^{k+1}}
\newcommand{\RT}{R_T^k}
\newcommand{\Rh}{R_h^k}
\newcommand{\VRT}{\bvec{R}_T^k}
\newcommand{\VRh}{\bvec{R}_h^k}

\newcommand{\PSigmaT}{\bvec{P}_{\bvec{\Sigma},T}^k}
\newcommand{\uGT}{\uvec{G}_T^k}
\newcommand{\uGh}{\uvec{G}_h^k}
\newcommand{\uRT}{\underline{R}_T^k}
\newcommand{\uRE}{\underline{R}_E^k}
\newcommand{\uRF}{\underline{R}_F^k}
\newcommand{\uRh}{\underline{R}_h^k}

\newcommand{\Err}{\mathcal{E}_h^k}
\newcommand{\dErotrot}{\widetilde{\mathcal{E}}_{\VROT\ROT,h}^k}
\newcommand{\Egrad}{\mathcal{E}_{\GRAD,h}^k}
\newcommand{\dEgrad}{\widetilde{\mathcal{E}}_{\GRAD,h}^k}

\begin{document}

\title{An arbitrary-order discrete rot-rot complex on polygonal meshes with application to a quad-rot problem}
\author[1]{Daniele A. Di Pietro}
\affil[1]{IMAG, Univ Montpellier, CNRS, Montpellier, France, \email{daniele.di-pietro@umontpellier.fr}}

\maketitle

\begin{abstract}
  In this work, following the discrete de Rham (DDR) approach, we develop a discrete counterpart of a two-dimensional de Rham complex with enhanced regularity.
  The proposed construction supports general polygonal meshes and arbitrary approximation orders.
  We establish exactness on a contractible domain for both the versions of the complex with and without boundary conditions and, for the former, prove a complete set of Poincar\'e-type inequalities.
  The discrete complex is then used to derive a novel discretisation method for a quad-rot problem which, unlike other schemes in the literature, does not require the forcing term to be prepared.
  We carry out complete stability and convergence analyses for the proposed scheme and provide numerical validation of the results.
\end{abstract}

%% \tableofcontents

%------------------------------------------------------------------------------%

\section{Introduction}

Denote by $\Omega\subset\Real^2$ a bounded connected polygonal set, which we assume contractible and of unit diameter for the sake of simplicity.
Following the discrete de Rham (DDR) approach of \cite{Di-Pietro.Droniou.ea:20,Di-Pietro.Droniou:21*1}, we develop a discrete counterpart of the following exact complex:
\begin{equation}\label{eq:continuous.complex}
  \begin{tikzcd}
    \Real
    \arrow[r,hook] & H^1(\Omega)
    \arrow{r}[above=2pt]{\GRAD} & \Hrotrot{\Omega}
    \arrow{r}[above=2pt]{\ROT} & H^1(\Omega)
    \arrow{r}[above=2pt]{0} & 0,
  \end{tikzcd}
\end{equation}
where $H^1(\Omega)$ denotes the space of scalar-valued functions that are square-integrable along with their gradient,
while $\Hrotrot{\Omega}\coloneq\left\{\bvec{v}\in\bvec{L}^2(\Omega;\Real^2)\st\text{$\ROT\bvec{v}\in L^2(\Omega)$ and $\VROT\ROT\bvec{v}\in\bvec{L}^2(\Omega;\Real^2)$}\right\}$; see \eqref{eq:2d.differential.operators} below for a definition of the scalar and vector two-dimensional rotors.
The complex \eqref{eq:continuous.complex} is relevant in the design of numerical schemes for quad-rot problems which arise in applications related to electromagnetism; see, e.g., \cite{Monk.Sun:12,Zheng.Hu.ea:11,Sun:16}.
A particularly delicate issue in the context of finite element approximations is the design of $\Hrotrot{\Omega}$-conforming spaces.
Two families of rot-rot-conforming finite elements that are among the first of this kind have been recently developed in \cite{Zhang.Wang.ea:19,Wang.Zhang.ea:22}.
In \cite{Hu.Zhang.ea:20}, these families, along with new ones, have been unified, fitted into complexes, and extended to the lowest-order case.

The above-cited constructions are based on conforming meshes composed of triangles or rectangles.
The possibility to support more general element shapes can be, however, a crucial advantage in the presence of singularities or complicated domain geometries, as it can be exploited to perform non-conforming mesh refinement or coarsening \cite{Antonietti.Giani.ea:13,Bassi.Botti.ea:12,Bassi.Botti.ea:14}.
General polygonal elements are supported by the virtual elements of \cite{Zhao.Zhang:21}, which hinge on standard differential operators and spaces of functions that are only partially computable.
In the present work, we follow instead a fully discrete approach that consists in replacing both the spaces and differential operators by discrete, entirely computable counterparts (see \cite{Beirao-da-Veiga.Dassi.ea:22} for an in-depth study of the relations between the virtual and fully discrete approaches).
The discrete complex proposed here is built starting from the two-dimensional DDR complex of \cite{Di-Pietro.Droniou:21*1} and has degrees of freedom (DOFs) that coincide with those of the most efficient complex in \cite{Zhao.Zhang:21} for the lowest-order version (corresponding, in the present context, to $k=0$), but differ at higher-orders.
Efficiency of the present complex can be enhanced using the serendipity techniques of \cite{Di-Pietro.Droniou:22}, resulting in most situations in a slighly lower DOF count with respect to \cite{Zhao.Zhang:21} as far as the discrete $H^1(\Omega)$ space at the head of the complex and the $\Hrotrot{\Omega}$ space are concerned.
The discrete $H^1(\Omega)$ space at the tail of the present complex, on the other hand, is smaller for $k\le 2$ or when elements with a sufficient number of sides are considered.
Another important difference is that the present construction is based on Koszul rather than orthogonal complements.
As noted in \cite{Di-Pietro.Droniou:21*1}, Koszul complements are usually easier to construct in practice and are hierarchical, a property which simplifies several arguments in the analysis. 

We present a full set of Poincar\'e inequalities for the version of the complex with boundary conditions.
Such inequalities lie at the core of the stability analysis of a novel numerical scheme for quad-rot problems based on our fully discrete complex.
We perform a complete study of this scheme, proving convergence in $h^{k+1}$ (with $h$ denoting, as usual, the meshsize) for the graph norm of the error when the complex of polynomial degree $k\ge 0$ is used as a starting point.
Convergence is numerically demonstrated to be in $h^{k+2}$ when only the $L^2$-like norms of the errors are considered.
Crucially, unlike other results in the literature \cite{Sun:16,Zhao.Zhang:21}, our scheme, and the corresponding analysis, do not require the forcing term to be prepared (i.e., either divergence-free or expressed in terms of a potential).
Working under the more realistic assumption that a potential for the forcing term is not available has deep implications in the analysis and requires, in particular, to prove an adjoint consistency result for rot-rot, which appears to be the first of this kind for polygonal methods.
Finally, it is also worth noting that, while we assume mesh elements to be contractible in order to ensure the exactness of the local complex, our analysis hinges on the results of \cite[Chapter 1]{Di-Pietro.Droniou:20} and therefore covers also non-star-shaped elements; see, in particular, Section 1.4 therein.
The theoretical results described above are corroborated by numerical evidence.

The rest of this work is organised as follows.
In Section \ref{sec:setting} we briefly outline the setting and recall the main notations.
Section \ref{sec:discrete.complex} contains the definition of the discrete complex and of its variants with boundary conditions and with serendipity-based DOFs reduction.
A detailed proof of the exactness for the complexes with and without boundary conditions is also provided.
In Section \ref{sec:poincare} we prove a complete set of Poincar\'e inequalities for the complex with boundary conditions, which is used as a starting point in Section \ref{sec:application} to design a stable and optimally convergent numerical scheme for a quad-rot problem.
 
%------------------------------------------------------------------------------%

\section{Setting}\label{sec:setting}

\subsection{Two-dimensional vector calculus operators}

Consider the real plane $\Real^2$ endowed with the Cartesian coordinate system $(x_1,x_2)$, and denote by $\partial_i$ the weak partial derivative with respect to the $i$th coordinate.
We need the following two-dimensional differential operators acting on smooth enough
scalar-valued fields $q$ or
vector-valued fields $\bvec{v}=\begin{pmatrix}v_1\\v_2\end{pmatrix}$:
\begin{equation}\label{eq:2d.differential.operators}
  \GRAD q\coloneq\begin{pmatrix}\partial_1 q\\ \partial_2 q\end{pmatrix},\quad
  \VROT q\coloneq\begin{pmatrix}\partial_2 q\\ -\partial_1 q\end{pmatrix}, \quad
  \ROT \bvec{v} \coloneq \partial_1 v_2 - \partial_2 v_1,\quad
  \DIV\bvec{v}\coloneq \partial_1 v_1 + \partial_2 v_2.
\end{equation}

\subsection{Mesh and notation for inequalities up to a constant}\label{sec:setting:mesh}

We denote by $\Mh = \Th\cup\Eh\cup\Vh$ a polygonal mesh of $\Omega$ in the usual sense of \cite{Di-Pietro.Droniou:20}, with $\Th$, $\Eh$, and $\Vh$ collecting, respectively, the elements, edges, and vertices and $h$ denoting the meshsize.
For all $Y\in\Mh$, we let $h_Y$ denote its diameter so that, in particular, $h=\max_{T\in\Th}h_T$.
$\Mh$ is assumed to belong to a refined mesh sequence with regularity parameter bounded away from zero.
We additionally assume that each element $T\in\Th$ is contractible, and denote by $\bvec{x}_T$ a point inside $T$ such that there exists a disk contained in $T$ centered in $\bvec{x}_T$ and of diameter comparable to $h_T$ uniformly in $h$.
The sets of edges and vertices of $T$ are denoted by $\ET$ and $\VT$, respectively.
By mesh regularity, the number of edges (and vertices) of mesh elements are bounded uniformly in $h$.
For each edge $E\in\Eh$, we denote by $\VE$ the set of vertices corresponding to its endpoints and fix an orientation by prescribing a unit tangent vector $\tangent_E$.
The corresponding unit normal vector $\normal_E$ is selected so that $(\tangent_E,\normal_E)$ forms a right-handed system of coordinates, and, for each $T\in\Th$ such that $E\in\ET$, we denote by $\omega_{TE}\in\{-1,+1\}$ the orientation of $E$ relative to $T$, defined so that $\omega_{TE}\normal_E$ points out of $T$.

From this point on, $a \lesssim b$ means $a\le Cb$ with $C$ only depending on $\Omega$, the mesh regularity parameter, and the polynomial degree $k$ of the discrete complex (see Section \ref{sec:discrete.complex}).
We also write $a\simeq b$ as a shorthand for ``$a\lesssim b$ and $b\lesssim a$''.

\subsection{Polynomial spaces}

Given $Y\in\Th\cup\Eh$ and an integer $m\ge 0$, we denote by $\Poly{m}(Y)$ the space spanned by the restriction to $Y$ of two-variate polynomials of total degree $\le m$,  with the additional convention that $\Poly{-1}(Y) \coloneq \{0\}$.
The symbol $\vPoly{m}(Y;\Real^2)$ denotes the set of vector-valued functions over $Y$ whose components are in $\Poly{m}(Y)$.
Finally we denote by $\Poly{m}(\Eh)$ the space of broken polynomials of total degree $\le m$ on $\Eh$ and, for all $T\in\Th$, by $\Poly{m}(\ET)$ its restriction to $\ET$.
Vector versions of this space are denoted in boldface and the codomain is specified.
We will need the following direct decompositions of $\vPoly{m}(T;\Real^2)$ (see, e.g., \cite[Corollary 7.4]{Arnold:18}):
\begin{equation}\label{eq:vPoly:decomposition}
  \vPoly{m}(T;\Real^2)
  = \Roly{m}(T) \oplus \cRoly{m}(T)
  = \Goly{m}(T) \oplus \cGoly{m}(T),,
\end{equation}
where
\[
\begin{alignedat}{4}
\Roly{m}(T)&\coloneq\VROT\Poly{m+1}(T),&\qquad&
\cRoly{m}(T)&\coloneq(\bvec{x} - \bvec{x}_T)\Poly{m-1}(T),
\\
\Goly{m}(T)&\coloneq\GRAD\Poly{m+1}(T),&\qquad&
\cGoly{m}(T)&\coloneq(\bvec{x} - \bvec{x}_T)^\top\Poly{m-1}(T),
\end{alignedat}
\]
where, for all $\bvec{v} = \begin{pmatrix}v_1\\v_2\end{pmatrix}\in\Real^2$, $\bvec{v}^\top \coloneq \begin{pmatrix}v_2\\-v_1\end{pmatrix}\in\Real^2$ is the vector obtained rotating $\bvec{v}$ by an angle of $-\frac\pi2$ radians.
Given a polynomial (sub)space $\mathcal{X}^m(Y)$ on $Y\in\Th\cup\Eh$, the corresponding $L^2$-orthogonal projector is denoted by $\pi_{\mathcal{X},Y}^m$.
Boldface fonts will be used when the elements of $\mathcal{X}^m(Y)$ are vector-valued, and we additionally denote by $\bvec{\pi}_{\cvec{X},T}^{\compl, m}$, $\cvec{X}\in\{\cvec{R}, \cvec{G}\}$, the $L^2$-orthogonal projector on $\cvec{X}^{\compl,m}(T)$.
The set of broken polynomials of total degree $\le m$ on the mesh is denoted by $\Poly{m}(\Th)$, and its vector version by $\vPoly{m}(\Th;\Real^2)$.
The $L^2$-orthogonal projector on $\Poly{m}(\Th)$ is denoted by $\lproj{k}{h}$.

%------------------------------------------------------------------------------%

\section{Discrete complex}\label{sec:discrete.complex}

Given a polygonal mesh $\Mh$ of $\Omega$ and a polynomial degree $k\ge 0$, the discrete version of \eqref{eq:continuous.complex} reads
\begin{equation}\label{eq:discrete.complex}
  \begin{tikzcd}
    \Real
    \arrow{r}[above=2pt]{\IV{h}} & \SV{h}
    \arrow{r}[above=2pt]{\uGh} & \SSigma{h}
    \arrow{r}[above=2pt]{\uRh} & \SW{h}
    \arrow{r}[above=2pt]{0} & 0,
  \end{tikzcd}
\end{equation}
where $\SV{h}$ denotes the standard two-dimensional counterpart of the $H^1(\Omega)$ space in the DDR construction,
$\SSigma{h}$ is a novel $\Hrotrot{\Omega}$-like space that constitutes the main novelty of this paper,
$\SW{h}$ is similar to $\SV{h}$ but with element polynomial components one degree higher,
and the operators $\IV{h}$, $\uGh$, and $\uRh$ are respectively defined by \eqref{eq:IV}, \eqref{eq:uGh}, and \eqref{eq:uRh} below.

\subsection{Discrete spaces}

The discrete spaces in \eqref{eq:discrete.complex} are defined as follows (see Remark \ref{rem:degree.SWh} below for further insight into the choice of the degree for element unknowns in $\SW{h}$):
\begin{equation}\label{eq:discrete.spaces}
  \begin{aligned}
    \SV{h} &\coloneq\Xgrad{k-1,k}{h},
    \\
    \SSigma{h} &\coloneq\Big\{
    \uvec{v}_h
    \begin{aligned}[t]
      &= \big(
      (\bvec{v}_{\cvec{R},T}, \bvec{v}_{\cvec{R},T}^\compl)_{T\in\Th},
      (v_E, C_{\bvec{v},E})_{E\in\Eh},
      (C_{\bvec{v},\nu})_{\nu\in\Vh}
      \big)\st
      \\
      &\text{%
        $\bvec{v}_{\cvec{R},T}\in\Roly{k-1}(T)$ and $\bvec{v}_{\cvec{R},T}^\compl\in\cRoly{k}(T)$ for all $T\in\Th$,
      }
      \\
      &\text{%
        $v_E\in\Poly{k}(E)$ and $C_{\bvec{v},E}\in\Poly{k-1}(E)$ for all $E\in\Eh$,
      }
      \\
      &\text{%
        $C_{\bvec{v},\nu}\in\Real$ for all $\nu\in\Vh$
      }\Big\},
    \end{aligned}
    \\
    \SW{h} &\coloneq\Xgrad{k,k}{h},
  \end{aligned}
\end{equation}
where, for $m\ge -1$,
\[
\Xgrad{m,k}{h}\coloneq\Big\{
\underline{q}_h
\begin{aligned}[t]
  &= \big( (q_T)_{T\in\Th}, (q_E)_{E\in\Eh}, (q_\nu)_{\nu\in\Vh}\big)\st
  \\
  &\text{%
    $q_T\in\Poly{m}(T)$ for all $T\in\Th$,
  }
  \\
  &\text{%
    $q_E\in\Poly{k-1}(F)$ for all $E\in\Eh$,
  }
  \\
  &\text{%
    $q_\nu\in\Real$ for all $\nu\in\Vh$
  }\Big\}.
\end{aligned}
\]
We denote the restrictions of the above spaces and of the operators mapping on them to $Y\in\Th\cup\Eh$ by replacing the subscript ``$h$'' with ``$Y$''.
Such restrictions are obtained collecting the components on $Y$ and its boundary.

\begin{remark}[Comparison between $\SSigma{h}$ and $\Xrot{h}$]\label{rem:comparison.Xcurl}
  Recall the definition of the standard two-dimensional DDR counterpart of the space $\Hrot{\Omega}\coloneq\left\{\bvec{v}\in\bvec{L}^2(\Omega;\Real^2)\st\ROT\bvec{v}\in L^2(\Omega)\right\}$:
  \[
  \Xrot{h}\coloneq\Big\{
  \uvec{v}_h
  \begin{aligned}[t]
    &= \big(
    (\bvec{v}_{\cvec{R},T}, \bvec{v}_{\cvec{R},T}^\compl)_{T\in\Th},
    (v_E)_{E\in\Eh}
    \big)\st
    \\
    &\text{%
      $\bvec{v}_{\cvec{R},T}\in\Roly{k-1}(T)$ and $\bvec{v}_{\cvec{R},T}^\compl\in\cRoly{k}(T)$ for all $T\in\Th$,
    }
    \\
    &\text{%
      $v_E\in\Poly{k}(E)$  for all $E\in\Eh$
    }\Big\}.
  \end{aligned}
  \]
  This space injects into $\SSigma{h}$ through the mapping
  \[
  \Injrot:\Xrot{h}\ni\uvec{v}_h\mapsto\big( (\bvec{v}_{\cvec{R},T}, \bvec{v}_{\cvec{R},T}^\compl)_{T\in\Th}, (v_E, 0)_{E\in\Eh}, (0)_{\nu\in\Vh}\big)\in\SSigma{h}.
  \]
  For future use, we also define the restriction
  \begin{equation}\label{eq:Resrot}
    \Resrot:\SSigma{h}\ni\uvec{v}_h\mapsto\big( (\bvec{v}_{\cvec{R},T},\bvec{v}_{\cvec{R},T}^\compl)_{T\in\Th}, (v_E)_{E\in\Eh}\big)\in\Xrot{h}.
  \end{equation}
\end{remark}

The meaning of the polynomial components in the spaces \eqref{eq:discrete.spaces} is provided by the interpolators.
Specifically, letting, for $\delta>0$,
\[
  \text{
    $V = W\coloneq H^{1+\delta}(\Omega)$ and $\bvec{\Sigma}\coloneq\big\{\bvec{v}\in\bvec{H}^{\nicefrac12+\delta}(\Omega)\st\ROT\bvec{v}\in H^{1+\delta}(\Omega)\big\}$,
  }
\]
we define $\IV{h}:V\to\SV{h}$ and $\ISigma{h}:\bvec{\Sigma}\to\SSigma{h}$ such that, for all $(q,\bvec{v},r)\in V\times\bvec{\Sigma}\times V$,
\begin{subequations}\label{eq:interpolators}
  \begin{align}\label{eq:IV}
    \IV{h} q
    &\coloneq\big(
    (\lproj{k-1}{T} q)_{T\in\Th}, (\lproj{k-1}{E} q)_{E\in\Eh}, (q(\bvec{x}_\nu))_{\nu\in\Vh}
    \big),
    \\ \label{eq:ISigma}
    \ISigma{h}\bvec{v}
    &\coloneq\big(
    (\Rproj{k-1}{T}\bvec{v},\cRproj{k}{T}\bvec{v})_{T\in\Th},
    (\lproj{k}{E}(\bvec{v}\cdot\tangent_E),\lproj{k-1}{E}(\ROT\bvec{v}))_{E\in\Eh},
    (\ROT\bvec{v}(\bvec{x}_\nu))_{\nu\in\Vh}
    \big),
    \\ \label{eq:IW}
    \IW{h} r
    &\coloneq\big(
    (\lproj{k}{T} r)_{T\in\Th}, (\lproj{k-1}{E} r)_{E\in\Eh}, (r(\bvec{x}_\nu))_{\nu\in\Vh}
    \big),
  \end{align}
\end{subequations}
where it is understood that local projectors are applied to restrictions or traces of their argument as needed.

\subsection{Discrete differential operators and potentials}

\subsubsection{Local gradients and scalar potentials}\label{sec:GT.PVT}

For all $E\in\Eh$, we let $\PVE:\SV{E}\to\Poly{k+1}(E)$ be such that, for all $\underline{q}_E$, $\PVE\underline{q}_E$ is the unique polynomial that satisfies $\PVE\underline{q}_E(\bvec{x}_\nu) = q_\nu$ for all $\nu\in\VE$ and $\lproj{k-1}{E}(\PVE\underline{q}_E) = q_E$.
The \emph{edge gradient} $\GE:\SV{E}\to\Poly{k}(E)$ is then simply defined setting $\GE\underline{q}_E \coloneq (\PVE\underline{q}_E)'$.

Let next $T\in\Th$ and $\underline{q}_T\in\SV{T}$.
The \emph{element gradient} $\GT:\SV{T}\to\vPoly{k}(T;\Real^2)$ and \emph{potential} $\PVT :\SV{T}\to\Poly{k+1}(T)$ are  obtained mimicking integration by parts formulas.
Specifically, given $\underline{q}_T\in\SV{T}$, $\GT\underline{q}_T$ and $\PVT\underline{q}_T$ are uniquely defined by the following conditions:
\begin{equation}\label{eq:GT}
  \int_T\GT\underline{q}_T\cdot\bvec{v}
  = -\int_T q_T~\DIV\bvec{v}
  + \sum_{E\in\ET}\omega_{TE}\int_E \PVE\underline{q}_E(\bvec{v}\cdot\normal_E)
  \qquad\forall\bvec{v}\in\vPoly{k}(T;\Real^2),
\end{equation}
and
\[
\int_T\PVT \underline{q}_T~\DIV\bvec{v}
= -\int_F\GT\underline{q}_T\cdot\bvec{v}
+ \sum_{E\in\ET}\omega_{TE}\int_E \PVE\underline{q}_E~(\bvec{v}\cdot\normal_E)
\qquad\forall\bvec{v}\in\cRoly{k+2}(T).
\]

\subsubsection{Local scalar rotor and vector potential}

Let again $T\in\Th$ and take $\uvec{v}_T\in\SSigma{T}$.
We define the \emph{element scalar rotor} $\RT:\SSigma{T}\to\Poly{k}(T)$ such that
\begin{equation}\label{eq:RT}
  \int_T\RT\uvec{v}_T~q
  = \int_T\bvec{v}_{\cvec{R},T}\cdot\VROT q
  - \sum_{E\in\ET}\omega_{TE}\int_E v_E~q
  \qquad\forall q\in\Poly{k}(T).
\end{equation}
As for the the element gradient, this definition mimics an integration by parts formula.
We have the following commutation property, which can be proved reasoning as in \cite[Proposition 4.3]{Di-Pietro.Droniou.ea:20}:
\begin{equation}\label{eq:RT:commutation}
  \RT\ISigma{T}\bvec{v} = \lproj{k}{T}(\ROT\bvec{v})
  \qquad\forall\bvec{v}\in\bvec{H}^1(T;\Real^2).
\end{equation}
At the global level, we let $\Rh:\Xrot{h}\to\Poly{k}(\Th)$ be such that, for all $\uvec{v}_h\in\Xrot{h}$,
\begin{equation}\label{eq:Rh}
  (\Rh\uvec{v}_h)_{|T}\coloneq\RT\uvec{v}_T
  \qquad\forall T\in\Th.
\end{equation}

The \emph{vector potential} $\PSigmaT:\SSigma{T}\to\vPoly{k}(T;\Real^2)$ is obtained, for a given $\uvec{v}_T\in\SSigma{T}$, taking $\PSigmaT\uvec{v}_T$ as the unique element of $\vPoly{k}(T;\Real^2)$ such that, for all $(q,\bvec{w})\in\Poly{k+1,0}(T)\times\cRoly{k}(T)$ (with $\Poly{k+1,0}(T)$ spanned by the functions of $\Poly{k+1}(T)$ with zero mean value over $T$),
\begin{equation}\label{eq:PSigmaT}
\int_T\PSigmaT\uvec{v}_T\cdot(\VROT q + \bvec{w})
= \int_T\RT\uvec{v}_T~q
+ \sum_{E\in\ET}\omega_{TE}\int_E v_E~q
+ \int_T\bvec{v}_{\cvec{R},T}^\compl\cdot\bvec{w}.
\end{equation}

\subsubsection{Global gradient and rotor}

The \emph{global gradient} $\uGh:\SV{h}\to\SSigma{h}$ and \emph{scalar rotor} $\uRh:\SSigma{h}\to\SW{h}$, acting between spaces of the discrete complex \eqref{eq:discrete.complex}, are defined setting, for all $(\underline{q}_h,\uvec{v}_h)\in\SV{h}\times\SSigma{h}$,
\begin{align}\label{eq:uGh}
  \uGh\underline{q}_h
  &\coloneq\big(
  (\Rproj{k-1}{T}\GT\underline{q}_T, \cRproj{k}{T}\GT\underline{q}_T)_{T\in\Th},
  (\GE\underline{q}_E, 0)_{E\in\Eh},
  (0)_{\nu\in\Vh}
  \big),
  \\ \label{eq:uRh}
  \uRh\uvec{v}_h
  &\coloneq\big(
  (\RT\uvec{v}_T)_{T\in\Th}, (C_{\bvec{v},E})_{E\in\Eh}, (C_{\bvec{v},\nu})_{\nu\in\Vh}
  \big).
\end{align}

\subsection{Exactness}

The goal of this section is to prove the following result:

\begin{theorem}[Exactness]\label{thm:exactness}
  The sequence \eqref{eq:discrete.complex} with spaces defined by \eqref{eq:discrete.spaces} and operators given by \eqref{eq:IV}, \eqref{eq:uGh}, and \eqref{eq:uRh} is an exact complex.
\end{theorem}

\begin{proof}
  \underline{1. \emph{Proof of $\Ker\uGh = \Image\Real$.}} This is an immediate consequence of the exactness of the global two-dimensional DDR complex (which follows from its local counterpart in \cite[Remark 10]{Di-Pietro.Droniou:21*1} proceeding along the lines of the proof of \cite[Theorem 3]{Di-Pietro.Droniou:21*1} since $\Omega$ is contractible).
  \smallskip\\
  \underline{2. \emph{Proof of $\Image\uGh=\Ker\uRh$.}}
  We start by proving that $\Image\uGh\subset\Ker\uRh$.
  To this purpose, it suffices to notice that, for all $\underline{q}_h\in\SV{h}$, $\RT\uGT\underline{q}_T = 0$ for all $T\in\Th$ by \cite[Proposition 2]{Di-Pietro.Droniou:21*1}, since $\RT$ only depends on the polynomial components shared by $\SSigma{T}$ and $\Xrot{T}$; see Remark \ref{rem:comparison.Xcurl}.
  The fact that the edge and vertex components of $\uRh\uGh\underline{q}_h$ are zero is an immediate consequence of the definition \eqref{eq:uGh} of the global gradient.

  Let us now prove the converse inclusion $\Ker\uRh\subset\Image\uGh$.
  Let $\uvec{v}_h\in\SSigma{h}$ be such that $\uRh\uvec{v}_h = \uvec{0}$, i.e.,
  \begin{equation}\label{eq:uRh=0}
    \text{%
      $\RT\uvec{v}_T = 0$ for all $T\in\Th$, $C_{\bvec{v},E} = 0$ for all $E\in\Eh$, and $C_{\bvec{v},\nu} = 0$ for all $\nu\in\Vh$.%
    }
  \end{equation}
  Recalling again Remark \ref{rem:comparison.Xcurl}, this means that there exists $\uvec{w}_h\in\Xrot{h}$ with vanishing standard discrete DDR scalar rotor such that $\uvec{v}_h = \Injrot\uvec{w}_h$.
  The exactness of the standard two-dimensional discrete de Rham complex then yields the existence of $\underline{q}_h\in\SV{h}$ such that $\uvec{w}_h = \big(
  (\Rproj{k-1}{T}\GT\underline{q}_T, \cRproj{k}{T}\GT\underline{q}_T)_{T\in\Th},
  (\GE\underline{q}_E)_{E\in\Eh}  
  \big)$.
  Recalling that $\uvec{v}_h = \Injrot\uvec{w}_h$ and using the definition \eqref{eq:uGh} of $\uGh$ proves that $\uvec{v}_h = \uGh\underline{q}_h$, which yields the desired inclusion.
  \smallskip\\
  \underline{3. \emph{Proof of the surjectivity of $\uRh$.}}
  By the exactness of the two-dimensional discrete de Rham complex, the global scalar rotor $\Rh$ defined by \eqref{eq:Rh}  is surjective in $\Poly{k}(\Th)$.
  To conclude the proof, it suffices to notice that the edge and vertex components of $\uRh$ span, respectively, $\Poly{k-1}(\Eh)$ and $\Real^{\card(\Vh)}$.
\end{proof}

\begin{remark}[Degree of the element components in $\SW{h}$]\label{rem:degree.SWh}
  Taking element components of degree $k$ in $\SW{k}$ is crucial to have the first relation in \eqref{eq:uRh=0}, which is in turn needed to leverage the exactness of the standard two-dimensional de Rham complex.
\end{remark}

\subsection{Cochain map property of the interpolators}

\begin{lemma}[Cochain map property of the interpolators]
  Denoting by $V$, $\bvec{\Sigma}$, and $W$ subspaces of the spaces appearing in the continuous complex \eqref{eq:continuous.complex} in which the interpolators \eqref{eq:interpolators} are well defined, the last two rows of the following picture form a commuting diagram:
  \begin{equation}\label{eq:commutation}
    \begin{tikzcd}[row sep=large]
      \Real
      \arrow[r,hook] & H^1(\Omega)
      \arrow{r}[above=2pt]{\GRAD} & \Hrotrot{\Omega}
      \arrow{r}[above=2pt]{\ROT} & H^1(\Omega)
      \arrow{r}[above=2pt]{0} & 0
      \\
      \Real
      \arrow[r,hook] & V\arrow{d}{\IV{h}}\arrow[u,hook]
      \arrow{r}[above=2pt]{\GRAD} & \bvec{\Sigma}\arrow{d}{\ISigma{h}}\arrow[u,hook]
      \arrow{r}[above=2pt]{\ROT} & W\arrow{d}{\IW{h}}\arrow[u,hook]
      \arrow{r}[above=2pt]{0} & 0
      \\
      \Real
      \arrow{r}[above=2pt]{\IV{h}} & \SV{h}
      \arrow{r}[above=2pt]{\uGh} & \SSigma{h}
      \arrow{r}[above=2pt]{\uRh} & \SW{h}
      \arrow{r}[above=2pt]{0} & 0.
    \end{tikzcd}
  \end{equation}
\end{lemma}

\begin{proof}
  The commutation property $\uGh\IV{h} q = \ISigma{h}(\GRAD q)$ for all $q\in V$ is a consequence of Point 1. in the proof of \cite[Lemma 4]{Di-Pietro.Droniou:21*1} along with the fact that the components of $\uGh\IV{h} q$ and $\ISigma{h}(\GRAD q)$ that are not in $\Xrot{h}$ are zero (the former by definition \eqref{eq:uGh} of the discrete gradient, the latter by definition \eqref{eq:ISigma} of $\ISigma{h}$ along with the fact that $\ROT\GRAD = 0$).
  \smallskip
  
  To prove the commutation property $\uRh\ISigma{h}\bvec{v} = \IW{h}(\ROT\bvec{v})$ for all $\bvec{v}\in\bvec{\Sigma}$, we write
  \[
  \begin{aligned}
    \uRh\ISigma{h}\bvec{v}
    &\stackrel{\text{\eqref{eq:uRh}, \eqref{eq:ISigma}}}{=}
    \big(
    (\RT\ISigma{h}\bvec{v})_{T\in\Th}, (\lproj{k-1}{E}(\ROT\bvec{v}))_{E\in\Eh}, (\ROT\bvec{v}(\bvec{x}_\nu))_{\nu\in\Vh}  
    \big)
    \\
    &\hspace{0.8em}\stackrel{\eqref{eq:RT:commutation}}{=}
    \big(
    (\lproj{k}{T}(\ROT\bvec{v}))_{T\in\Th}, (\lproj{k-1}{E}(\ROT\bvec{v}))_{E\in\Eh}, (\ROT\bvec{v}(\bvec{x}_\nu))_{\nu\in\Vh}  
    \big)
    \stackrel{\eqref{eq:IW}}{=}
    \IW{h}(\ROT\bvec{v}).\qedhere
  \end{aligned}
  \]
\end{proof}

\subsection{Discrete complex with boundary conditions}

In the application of Section \ref{sec:application} below, we will use a variant of the complex \eqref{eq:discrete.complex} with boundary conditions.
To this purpose, denoting by $\Ehb$ and $\Vhb$ the sets of boundary edges and vertices, we let
\begin{align}\notag
  \SVZ{h}
  &\coloneq\big\{
  \underline{q}_h\in\SV{h}\st
  \text{$q_E = 0$ for all $E\in\Ehb$ and $q_\nu = 0$ for all $\nu\in\Vhb$}
  \big\},
  \\ \label{eq:SSigmaZ}
  \SSigmaZ{h}
  &\coloneq\big\{
  \uvec{v}_h\in\SSigma{h}\st
  \text{$v_E = C_{\bvec{v},E} = 0$ for all $E\in\Ehb$ and $C_{\bvec{v},\nu} = 0$ for all $\nu\in\Vhb$}
  \big\},
  \\\notag
  \SWZ{h}
  &\coloneq\big\{
  \underline{r}_h\in\SV{h}\st
  \text{$r_E = 0$ for all $E\in\Ehb$ and $r_\nu = 0$ for all $\nu\in\Vhb$}
  \big\}.
\end{align}

\begin{theorem}[Exactness of the complex with boundary conditions]\label{thm:exactness.bc}
  The following sequence defines an exact complex:
  \begin{equation}\label{eq:discrete.complex:bc}
    \begin{tikzcd}
      0
      \arrow{r}[above=2pt]{\IV{h}} & \SVZ{h}
      \arrow{r}[above=2pt]{\uGh} & \SSigmaZ{h}
      \arrow{r}[above=2pt]{\uRh} & \SWZ{h}
      \arrow{r}[above=2pt]{0} & 0.
    \end{tikzcd}
  \end{equation}
\end{theorem}

\begin{proof}
  \underline{1. \emph{Complex property.}}
  Let's start to prove that \eqref{eq:discrete.complex:bc} is a complex.
  Clearly, $\IV{h}0 = \underline{0}\in\SVZ{h}$.
  Next, for all $\underline{q}_h\in\SVZ{h}$, it holds by definition that $\GE\underline{q}_E = 0$ for all $E\in\Ehb$ which, combined with the definition \eqref{eq:uGh} of the global discrete gradient, proves that $\uGh\underline{q}_h\in\SSigmaZ{h}$.
  Finally, for all $\uvec{v}_h\in\SSigmaZ{h}$, by definition \eqref{eq:SSigmaZ} of this space, the components on boundary edges and vertices of $\uRh\uvec{v}_h$ defined by \eqref{eq:uRh} vanish, proving that $\uRh\uvec{v}_h\in\SWZ{h}$.
  \medskip\\ 
  \underline{2. \emph{Exactness.}}
  We next prove that the complex \eqref{eq:discrete.complex:bc} is exact.
  Let $\underline{q}_h\in\SVZ{h}$ be such that $\uGh\underline{q}_h = \uvec{0}$.
  By exactness of the complex \eqref{eq:discrete.complex} without boundary conditions, there is $C\in\Real$ such that $\underline{q}_h = \IV{h} C$, and the only possibility compatible with the fact that the boundary components of $\underline{q}_h$ vanish is $C = 0$.
  \smallskip
  
  Let now $\uvec{v}_h\in\SSigmaZ{h}$ be such that $\uRh\uvec{v}_h = \underline{0}$.
  By exactness of the discrete complex without boundary conditions, there is $\underline{q}_h\in\SV{h}$ (defined up to an element of $\IV{h}\Real$) such that $\uvec{v}_h = \uGh\underline{q}_h$.
  Since $\uvec{v}_h\in\SSigmaZ{h}$, the boundary condition implies $\GE\underline{q}_E = 0$ for all $E\in\Ehb$ so that, by single-valuedness of vertex unknowns, there is $C\in\Real$ such that $q_\nu = C$ for all $\nu\in\Vhb$ and $q_E = \lproj{k-1}{E} C$ for all $E\in\Ehb$.
  Upon the substitution $\underline{q}_h - \IV{h} C$, we have thus found $\underline{q}_h\in\SVZ{h}$ such that $\uvec{v}_h = \uGh\underline{q}_h$.
  \smallskip
  
  To conclude, it only remains to prove that $\uRh:\SSigmaZ{h}\to\SWZ{h}$ is surjective.
  Let $\underline{r}_h\in\SWZ{h}$. By surjectivity of $\uRh:\SSigma{h}\to\SW{h}$, there exists $\uvec{v}_h\in\SSigma{h}$ (defined up to an element of $\uGh\SV{h}$) such that $\underline{r}_h = \uRh\uvec{v}_h$.
  By definition \eqref{eq:uRh} of the global scalar rotor, it holds $C_{\bvec{v},E} = 0$ for all $E\in\Ehb$ and $C_{\bvec{v},\nu} = 0$ for all $\nu\in\Vhb$.
  It only remains to make sure that we can take $v_E = 0$ for all $E\in\Ehb$.
  To this purpose, we start by noticing that the zero rotor condition implies
  \[
  0 = \sum_{T\in\Th}\int_T\RT\uvec{v}_h
  \stackrel{\eqref{eq:RT}}{=} -\sum_{T\in\Th}\sum_{E\in\ET}\omega_{TE}\int_E v_E
  = -\sum_{E\in\Ehb}\omega_{T_EE}\int_E v_E,
  \]
  where the conclusion follows observing that, by definition, $\omega_{T_1E} + \omega_{T_2E} = 0$ for all $E\in\Eh\setminus\Ehb$ shared by the mesh elements $T_1$ and $T_2$ (for all $E\in\Ehb$, $T_E\in\Th$ denotes the unique mesh element such that $E\in\edges{T_E}$).
  This condition implies the existence of $q_{\partial\Omega}\in\Poly{k+1}(\Ehb)\cap C^0(\partial\Omega)$ (with $\Poly{k+1}(\Ehb)$ spanned by broken polynomials of total degree $\le k+1$ on $\Ehb$) such that $\omega_{T_EE}v_E = (q_{\partial\Omega})_{|E}'$ for all $E\in\Ehb$.
  Substituting $\uvec{v}_h\gets\uvec{v}_h - \uGh\underline{q}_h$ with $\underline{q}_h$ lifting of $q_{\partial\Omega}$ such that $q_\nu = q_{\partial\Omega}(\bvec{x}_\nu)$ for all $\nu\in\Vhb$ and $q_E = \lproj{k-1}{E}(q_{\partial\Omega})_{|E}$ for all $E\in\Ehb$, we have thus found $\uvec{v}_h\in\SSigmaZ{h}$ such that $\underline{r}_h = \uRh\uvec{v}_h$ and concluded the proof.
\end{proof}

\subsection{Serendipity complex and comparison}\label{sec:discrete.complex:serendipity}

A leaner version of the discrete complex can be obtained using serendipity techniques as in \cite{Di-Pietro.Droniou:22}.
Specifically, for all $T\in\Th$, denote by $\eta_T\ge 2$ the number of not pairwise aligned edges selected to satisfy Assumption 11 therein, and set
\[
\ell_T \coloneq \max(k + 1 - \eta_T, -1)\le k-1.
\]
Then, the spaces $\SV{h}$ and $\SSigma{h}$ can be replaced with the following serendipity versions:
\[
\SVser{k}\coloneq\Xgrad{\ell_T,k}{h},\qquad
\SSigmaser{h} \coloneq\Big\{
\uvec{v}_h
\begin{aligned}[t]
  &= \big(
  (\bvec{v}_{\cvec{R},T}, \bvec{v}_{\cvec{R},T}^\compl)_{T\in\Th},
  (v_E, C_{\bvec{v},E})_{E\in\Eh},
  (C_{\bvec{v},\nu})_{\nu\in\Vh}
  \big)\st
  \\
  &\text{%
    $\bvec{v}_{\cvec{R},T}\in\Roly{k-1}(T)$ and $\bvec{v}_{\cvec{R},T}^\compl\in\cRoly{\ell_T+1}(T)$ for all $T\in\Th$,
  }
  \\
  &\text{%
    $v_E\in\Poly{k}(E)$ and $C_{\bvec{v},E}\in\Poly{k-1}(E)$ for all $E\in\Eh$,
  }
  \\
  &\text{%
    $C_{\bvec{v},\nu}\in\Real$ for all $\nu\in\Vh$
  }\Big\}.
\end{aligned}
\]
Introducing extension and reduction operators similar to those described in \cite[Section 5.3]{Di-Pietro.Droniou:22}, and using them to define interpolators and discrete gradient and rotor operators mapping on the serendipity spaces, one obtains again an exact complex.
For $k\ge 1$, and provided $\ell_T>2$ for at least some $T\in\Th$, such complex has fewer unknowns than the one in \eqref{eq:discrete.complex} owing to the reduction in the degree of element polynomial components.
The DOF count for both the full and serendipity DDR complexes is provided in Table \ref{tab:dofs.reduction}, where we also include a comparison with the most efficient construction in \cite{Zhao.Zhang:21} (corresponding to the choice $r=k-1$ therein).
Notice that a different convention is used here, so the polynomial degree with respect to the above reference is shifted by $-2$.

\begin{table}\centering
\begin{tabular}{c|ccccc}
\toprule
Discrete space & $k=0$ & $k=1$ & $k=2$ & $k=3$ & $k=4$ \\
\midrule
\multicolumn{ 6 }{c}{ Triangle, $\eta_T =  3 $ } \\
\midrule
$H^1(T)$ (head)
 & \cellcolor{black!10}{3 \textbullet{} 3 \textbullet{} 6} & \cellcolor{black!20}{7 \textbullet{} 6 \textbullet{} 10} & \cellcolor{black!20}{12 \textbullet{} 10 \textbullet{} 15} & \cellcolor{black!20}{18 \textbullet{} 15 \textbullet{} 21} & \cellcolor{black!20}{25 \textbullet{} 21 \textbullet{} 28} \\
$\Hrotrot{T}$
& \cellcolor{black!10}{6 \textbullet{} 6 \textbullet{} 6}& 15 \textbullet{} 14 \textbullet{} 13& \cellcolor{black!10}{26 \textbullet{} 23 \textbullet{} 23}& \cellcolor{black!20}{39 \textbullet{} 34 \textbullet{} 35}& \cellcolor{black!20}{54 \textbullet{} 47 \textbullet{} 49} \\
$H^1(T)$ (tail)
& \cellcolor{black!10}{4 \textbullet{} 4 \textbullet{} 6}& \cellcolor{black!10}{9 \textbullet{} 9 \textbullet{} 10}& \cellcolor{black!10}{15 \textbullet{} 15 \textbullet{} 15}& 22 \textbullet{} 22 \textbullet{} 21& 30 \textbullet{} 30 \textbullet{} 28 \\
\midrule
\multicolumn{ 6 }{c}{ Quadrangle, $\eta_T =  4 $ } \\
\midrule
$H^1(T)$ (head)
 & \cellcolor{black!10}{4 \textbullet{} 4 \textbullet{} 8} & \cellcolor{black!20}{9 \textbullet{} 8 \textbullet{} 13} & \cellcolor{black!20}{15 \textbullet{} 12 \textbullet{} 19} & \cellcolor{black!20}{22 \textbullet{} 17 \textbullet{} 26} & \cellcolor{black!20}{30 \textbullet{} 23 \textbullet{} 34} \\
$\Hrotrot{T}$
& \cellcolor{black!10}{8 \textbullet{} 8 \textbullet{} 8}& 19 \textbullet{} 18 \textbullet{} 17& \cellcolor{black!10}{32 \textbullet{} 29 \textbullet{} 29}& \cellcolor{black!20}{47 \textbullet{} 41 \textbullet{} 43}& \cellcolor{black!20}{64 \textbullet{} 55 \textbullet{} 59} \\
$H^1(T)$ (tail)
& \cellcolor{black!10}{5 \textbullet{} 5 \textbullet{} 8}& \cellcolor{black!10}{11 \textbullet{} 11 \textbullet{} 13}& \cellcolor{black!10}{18 \textbullet{} 18 \textbullet{} 19}& \cellcolor{black!10}{26 \textbullet{} 26 \textbullet{} 26}& 35 \textbullet{} 35 \textbullet{} 34 \\
\midrule
\multicolumn{ 6 }{c}{ Pentagon, $\eta_T =  5 $ } \\
\midrule
$H^1(T)$ (head)
 & \cellcolor{black!10}{5 \textbullet{} 5 \textbullet{} 10} & \cellcolor{black!20}{11 \textbullet{} 10 \textbullet{} 16} & \cellcolor{black!20}{18 \textbullet{} 15 \textbullet{} 23} & \cellcolor{black!20}{26 \textbullet{} 20 \textbullet{} 31} & \cellcolor{black!20}{35 \textbullet{} 26 \textbullet{} 40} \\
$\Hrotrot{T}$
& \cellcolor{black!10}{10 \textbullet{} 10 \textbullet{} 10}& 23 \textbullet{} 22 \textbullet{} 21& \cellcolor{black!10}{38 \textbullet{} 35 \textbullet{} 35}& \cellcolor{black!20}{55 \textbullet{} 49 \textbullet{} 51}& \cellcolor{black!20}{74 \textbullet{} 64 \textbullet{} 69} \\
$H^1(T)$ (tail)
& \cellcolor{black!10}{6 \textbullet{} 6 \textbullet{} 10}& \cellcolor{black!10}{13 \textbullet{} 13 \textbullet{} 16}& \cellcolor{black!10}{21 \textbullet{} 21 \textbullet{} 23}& \cellcolor{black!10}{30 \textbullet{} 30 \textbullet{} 31}& \cellcolor{black!10}{40 \textbullet{} 40 \textbullet{} 40} \\
\midrule
\multicolumn{ 6 }{c}{ Hexagon, $\eta_T =  6 $ } \\
\midrule
$H^1(T)$ (head)
 & \cellcolor{black!10}{6 \textbullet{} 6 \textbullet{} 12} & \cellcolor{black!20}{13 \textbullet{} 12 \textbullet{} 19} & \cellcolor{black!20}{21 \textbullet{} 18 \textbullet{} 27} & \cellcolor{black!20}{30 \textbullet{} 24 \textbullet{} 36} & \cellcolor{black!20}{40 \textbullet{} 30 \textbullet{} 46} \\
$\Hrotrot{T}$
& \cellcolor{black!10}{12 \textbullet{} 12 \textbullet{} 12}& 27 \textbullet{} 26 \textbullet{} 25& \cellcolor{black!10}{44 \textbullet{} 41 \textbullet{} 41}& \cellcolor{black!20}{63 \textbullet{} 57 \textbullet{} 59}& \cellcolor{black!20}{84 \textbullet{} 74 \textbullet{} 79} \\
$H^1(T)$ (tail)
& \cellcolor{black!10}{7 \textbullet{} 7 \textbullet{} 12}& \cellcolor{black!10}{15 \textbullet{} 15 \textbullet{} 19}& \cellcolor{black!10}{24 \textbullet{} 24 \textbullet{} 27}& \cellcolor{black!10}{34 \textbullet{} 34 \textbullet{} 36}& \cellcolor{black!10}{45 \textbullet{} 45 \textbullet{} 46} \\
\bottomrule
\end{tabular}
\caption{Number of DOFs for the full \textbullet{} serendipity \textbullet{} virtual \cite{Zhao.Zhang:21} (with $k$ shifted by $-2$) discrete counterparts of the spaces $H^1(T)$ (head of the complex), $\Hrotrot{T}$, and $H^1(T)$ (tail of the complex) on a triangle, quadrangle, pentagon, and hexagon element $T$ for polynomial degrees $k$ ranging from 0 to 4. The relative DOFs reduction is in parenthesis. The parameter $\eta_T$ is defined in Section \ref{sec:discrete.complex:serendipity}. The cases in which the serendipity DDR space is smaller (resp. smaller or equal) than the other two are highlighted in dark gray (resp. light gray). \label{tab:dofs.reduction}}
\end{table}

%------------------------------------------------------------------------------%

\section{Poincar\'e inequalities on the discrete complex with boundary conditions}\label{sec:poincare}

The stability of the scheme for the quad-rot problem considered in Section \ref{sec:application} below hinges on Poincar\'e inequalities on the discrete gradient and rot-rot operators proved in this section.

\subsection{Norms}\label{sec:poincare:norms}

We equip the spaces $\SV{h}$, $\Xrot{h}$, $\SSigma{h}$, and $\SW{h}$ with the standard $L^2$-like DDR \emph{component norms} such that, for all $(\underline{q}_h,\uvec{w}_h,\uvec{v}_h,\underline{r}_h)\in\SV{h}\times\Xrot{h}\times\SSigma{h}\times\SW{h}$,
\[
\begin{alignedat}{3}
  \tnorm[V,h]{\underline{q}_h}^2&\coloneq\sum_{T\in\Th}\tnorm[V,T]{\underline{q}_T}^2,
  &\qquad
  \tnorm[\ROT,h]{\uvec{w}_h}^2&\coloneq\sum_{T\in\Th}\tnorm[\ROT,T]{\uvec{w}_T}^2,
  \\
  \tnorm[\bvec{\Sigma},h]{\uvec{v}_h}^2&\coloneq\sum_{T\in\Th}\tnorm[\bvec{\Sigma},T]{\uvec{v}_T}^2,
  &\qquad
  \tnorm[V,h]{\underline{r}_h}^2&\coloneq\sum_{T\in\Th}\tnorm[V,T]{\underline{r}_T}^2,
\end{alignedat}
\]
where, for all $T\in\Th$,
\begin{align}\notag
  \tnorm[V,T]{\underline{q}_T}^2
  &\coloneq
  \norm[L^2(T)]{q_T}^2
  + h_T\sum_{E\in\ET}\norm[L^2(E)]{q_E}^2
  + h_T^2\sum_{\nu\in\VT}|q_\nu|^2,
  \\\label{eq:tnorm.rot.T}
  \tnorm[\ROT,T]{\uvec{w}_T}^2
  &\coloneq
  \norm[\bvec{L}^2(T;\Real^2)]{\bvec{w}_{\cvec{R},T}}^2
  + \norm[\bvec{L}^2(T;\Real^2)]{\bvec{w}_{\cvec{R},T}^\compl}^2
  + h_T\sum_{E\in\ET}\norm[L^2(E)]{w_E}^2,
  \\ \label{eq:tnorm.Sigma.T}
  \tnorm[\bvec{\Sigma},T]{\uvec{v}_T}^2
  &\coloneq
  \tnorm[\ROT,T]{\ResrotT\uvec{v}_T}^2
  + h_T^3\sum_{E\in\ET}\norm[L^2(E)]{C_{\bvec{v},E}}^2
  +  h_T^4 \sum_{\nu\in\VT}|C_{\bvec{v},\nu}|^2,
  \\ \label{eq:tnorm.W.h}
  \tnorm[W,T]{\underline{r}_T}^2
  &\coloneq
  \norm[L^2(T)]{r_T}^2
  + h_T\sum_{E\in\ET}\norm[L^2(E)]{r_E}^2
  + h_T^2\sum_{\nu\in\VT}|r_\nu|^2,
\end{align} 
with $\ResrotT:\SSigma{T}\ni\uvec{v}_T\to\big(\bvec{v}_{\cvec{R},T},\bvec{v}_{\cvec{R},T}^\compl,(v_E)_{E\in\ET}\big)\in\Xrot{T}$ local counterpart of the restriction \eqref{eq:Resrot}.
Above, the factors involving powers of $h_T$ ensure that all the terms have the same scaling.

On $\SSigma{h}$, we will also need the \emph{operator norm} induced by the following discrete $L^2$-product:
For all $(\uvec{w}_h,\uvec{v}_h)\in\SSigma{h}\times\SSigma{h}$,
\begin{equation}\label{eq:inner.prod:Sigma.h}
  (\uvec{w}_h,\uvec{v}_h)_{\bvec{\Sigma},h}
  \coloneq
  \sum_{T\in\Th}(\uvec{w}_T,\uvec{v}_T)_{\bvec{\Sigma},T},
\end{equation}
where, for all $T\in\Th$,
\begin{equation}\label{eq:inner.prod:Sigma.T}
  \begin{aligned}
    (\uvec{w}_T,\uvec{v}_T)_{\bvec{\Sigma},T}
    &\coloneq
    \int_T\PSigmaT\uvec{w}_T\cdot\PSigmaT\uvec{v}_T
    +  h_T\sum_{E\in\ET}\int_E(\PSigmaT\uvec{w}_T\cdot\tangent_E - w_E)~(\PSigmaT\uvec{v}_T\cdot\tangent_E - v_E)
    \\
    &\quad
    + h_T^3\sum_{E\in\ET}\int_E\lproj{k-1}{E}(\RT\uvec{w}_T - C_{\bvec{w},E})~\lproj{k-1}{E}(\RT\uvec{v}_T - C_{\bvec{v},E})
    \\
    &\quad
    + h_T^4\sum_{\nu\in\VT}(\RT\uvec{w}_T(\bvec{x}_\nu) - C_{\bvec{w},\nu})~(\RT\uvec{v}_T(\bvec{x}_\nu) - C_{\bvec{v},\nu}).
  \end{aligned}
\end{equation}
The induced local and global norms are, for all $\bullet\in\Th\cup\{h\}$ and all $\uvec{v}_\bullet\in\SSigma{\bullet}$,
\begin{equation}\label{eq:norm.Sigma}
  \norm[\bvec{\Sigma},\bullet]{\uvec{v}_\bullet}\coloneq
  (\uvec{v}_\bullet,\uvec{v}_\bullet)_{\bvec{\Sigma},\bullet}^{\nicefrac12}.
\end{equation}

\begin{lemma}[Equivalence of norms in $\SSigma{T}$ and $\SSigma{h}$]
  It holds, for $\bullet\in\Th\cup\{h\}$,
  \begin{equation}\label{eq:norm.equivalence:SSigma.h}
    \tnorm[\bvec{\Sigma},\bullet]{\uvec{v}_\bullet}
    \simeq\norm[\bvec{\Sigma},\bullet]{\uvec{v}_\bullet}
    \qquad\forall\uvec{v}_\bullet\in\SSigma{\bullet}.
  \end{equation}
\end{lemma}

\begin{proof}
  It suffices to prove the result for a generic $\bullet = T\in\Th$, as the result for $\bullet = h$ follows squaring the latter, summing over $T\in\Th$, and taking the square root of the resulting equivalence.
  The two-dimensional counterpart of the norm equivalence in \cite[Lemma~5]{Di-Pietro.Droniou:21*1} gives, for all $\uvec{v}_T\in\SSigma{T}$,
  \begin{equation}\label{eq:equivalence:tnorm.norm.rot}
    \tnorm[\ROT,T]{\ResrotT\uvec{v}_T}^2\simeq
    \norm[\bvec{L}^2(T;\Real^2)]{\PSigmaT\uvec{v}_T}^2
    + h_T\sum_{E\in\ET}\norm[L^2(E)]{\PSigmaT\uvec{v}_T\cdot\tangent_E - v_E}^2
    \eqcolon\norm[\ROT,T]{\Resrot\uvec{v}_T}^2.
  \end{equation}
  Moreover, for all $T\in\Th$, taking $q = \RT\uvec{v}_T$ in the definition \eqref{eq:RT} of the element scalar rotor, using triangle, inverse, and trace inequalities in the right-hand side, and simplifying, we obtain
  \begin{equation}\label{eq:RT:boundedness}
    \norm[L^2(T)]{\RT\uvec{v}_T}
    \lesssim h_T^{-1}\tnorm[\ROT,T]{\Resrot\uvec{v}_T}
    \stackrel{\eqref{eq:equivalence:tnorm.norm.rot}}{\lesssim}
    h_T^{-1}\norm[\ROT,T]{\Resrot\uvec{v}_T}
    \qquad\forall\uvec{v}_T\in\SSigma{T}.
  \end{equation}
  \medskip
  \underline{(i) \emph{Proof of $\tnorm[\bvec{\Sigma},T]{\uvec{v}_T}\lesssim\norm[\bvec{\Sigma},T]{\uvec{v}_T}$.}}
  For all $E\in\ET$, inserting $\pm\RT\uvec{v}_T$ into the norm, and using a triangle inequality followed by the continuity of $\lproj{k-1}{E}$ and a discrete trace inequality, we obtain
  \begin{equation}\label{eq:est.CvE}
    \begin{aligned}
      \norm[L^2(E)]{C_{\bvec{v},E}}
      &\le\norm[L^2(E)]{\lproj{k-1}{E}\RT\uvec{v}_T - C_{\bvec{v},E}}
      + h_T^{-\nicefrac12}\norm[L^2(T)]{\RT\uvec{v}_T}
      \\
      &\!\stackrel{\eqref{eq:RT:boundedness}}{\lesssim}
      \norm[L^2(E)]{\lproj{k-1}{E}\RT\uvec{v}_T - C_{\bvec{v},E}}
      + h_T^{-\nicefrac32}\norm[\ROT,T]{\Resrot\uvec{v}_T}.
    \end{aligned}
  \end{equation}
  Proceeding similarly, for all $\nu\in\VT$ we get
  \begin{equation}\label{eq:est.CvV}
    |C_{\bvec{v},\nu}|    
    \lesssim
    |\RT\uvec{v}_T(\bvec{x}_\nu) - C_{\bvec{v},\nu}|
    + h_T^{-2}\norm[\ROT,T]{\Resrot\uvec{v}_T}.
  \end{equation}
  Plugging \eqref{eq:est.CvE} and \eqref{eq:est.CvV} into the definition \eqref{eq:tnorm.Sigma.T} of $\tnorm[\bvec{\Sigma},T]{\uvec{v}_T}$, using the fact that the number of edges and vertices of $T$ is bounded uniformly in $h$ by mesh regularity, and recalling \eqref{eq:equivalence:tnorm.norm.rot}, the conclusion follows.
  \medskip\\
  \underline{(ii) \emph{Proof of $\norm[\bvec{\Sigma},T]{\uvec{v}_T}\lesssim\tnorm[\bvec{\Sigma},T]{\uvec{v}_T}$.}}
  For all $E\in\ET$, using a triangle inequality followed by the continuity of $\lproj{k-1}{T}$, a discrete trace inequality, and \eqref{eq:RT:boundedness}, we obtain
  \begin{equation}\label{eq:pen.E}    
    \norm[L^2(E)]{\lproj{k-1}{T}\RT\uvec{v}_T - C_{\bvec{v},E}}
    \lesssim\norm[L^2(E)]{C_{\bvec{v},E}}
    + h_T^{-\nicefrac32}\tnorm[\ROT,T]{\Resrot\uvec{v}_T}.
  \end{equation}
  In a similar way, for all $\nu\in\VT$ we get
  \begin{equation}\label{eq:pen.V}
    |\RT\uvec{v}_T(\bvec{x}_\nu) - C_{\bvec{v},\nu}|
    \lesssim
    |C_{\bvec{v},\nu}|
    + h_T^{-2}\tnorm[\ROT,T]{\Resrot\uvec{v}_T}.
  \end{equation}
    Plugging \eqref{eq:pen.E} and \eqref{eq:pen.V} into the definition of $\norm[\bvec{\Sigma},T]{\uvec{v}_T}$, using the fact that the number of edges and vertices of $T$ is bounded uniformly in $h$ by mesh regularity, and recalling \eqref{eq:equivalence:tnorm.norm.rot}, the conclusion follows.
\end{proof}
  
\subsection{Poincar\'e inequality in $\SVZ{h}$}

\begin{lemma}[Poincar\'e inequality in $\SVZ{h}$]\label{lem:poincare:VZh}
  For all $\underline{q}_h\in\SVZ{h}$, it holds
  \begin{equation}\label{eq:poincare:VZh}
    \tnorm[V,h]{\underline{q}_h}
    \lesssim\tnorm[\bvec{\Sigma},h]{\uGh\underline{q}_h}.
  \end{equation}
\end{lemma}

\begin{proof}
  Let $\hat{\underline{q}}_h\coloneq\big( (\PVT \underline{q}_T)_{T\in\Th}, (\PVE\underline{q}_E)_{E\in\Eh}\big)$ and notice that its boundary components are zero since $\underline{q}_h\in\SVZ{h}$.
  By the discrete Poincaré inequality for HHO methods \cite[Lemma 2.15]{Di-Pietro.Droniou:20}, it holds
  \begin{equation}\label{eq:poincare:VZh:1}
    \begin{aligned}
      \sum_{T\in\Th}\norm[L^2(T)]{\PVT\underline{q}_T}^2
      &\lesssim\sum_{T\in\Th}\left(
      \norm[\bvec{L}^2(T;\Real^2)]{\GRAD\PVT\underline{q}_T}^2
      + h_T^{-1}\!\sum_{E\in\ET}\norm[L^2(E)]{\PVT\underline{q}_T - \PVE\underline{q}_E}^2
      \right)
      \\
      &\lesssim\tnorm[\ROT,h]{\Resrot\uGh\underline{q}_h}^2
      \le\tnorm[\bvec{\Sigma},h]{\uGh\underline{q}_h}^2,
    \end{aligned}
  \end{equation}
  where we have used \cite[Eq.~(5.5)]{Di-Pietro.Droniou:21*1} to pass to the second line
  and observed that $\tnorm[\ROT,h]{\Resrot\uvec{v}_h}\le\tnorm[\bvec{\Sigma},h]{\uvec{v}_h}$ (cf. \eqref{eq:tnorm.Sigma.T}) to conclude.
  Noticing that $h_T\le 1\le h_T^{-1}$ (since the diameter of $\Omega$ is equal to 1 by assumption and $h_T$ is bounded by the latter quantity), and using the \eqref{eq:poincare:VZh:1}, we moreover have
  \begin{equation}\label{eq:poincare:VZh:2}
    \sum_{T\in\Th} h_T\sum_{E\in\ET}\norm[L^2(E)]{\PVT\underline{q}_T - \PVE\underline{q}_E}^2
    %% \le\sum_{T\in\Th} h_T^{-1}\sum_{E\in\ET}\norm[L^2(E)]{\PVT\underline{q}_T - \PVE\underline{q}_E}^2
    \lesssim\tnorm[\bvec{\Sigma},h]{\uGh\underline{q}_h}^2.
  \end{equation}
  Summing \eqref{eq:poincare:VZh:1} and \eqref{eq:poincare:VZh:2}, and using a norm equivalence which is the two-dimensional counterpart of \cite[Lemma 5]{Di-Pietro.Droniou:21*1} to bound the term in the left-hand side of the resulting inequality from below, we finally get  
  \[
  \tnorm[V,h]{\underline{q}_h}^2
  \lesssim
  \sum_{T\in\Th}\left(
  \norm[L^2(T)]{\PVT\underline{q}_T}^2
  + h_T\sum_{E\in\ET}\norm[L^2(E)]{\PVT\underline{q}_T - \PVE\underline{q}_E}^2
  \right)
  \lesssim\tnorm[\bvec{\Sigma},h]{\uGh\underline{q}_h}^2.\qedhere
  \]
\end{proof}

\subsection{Poincar\'e inequality in $\XrotZ{h}$}

As an intermediate step to prove a Poincar\'e inequality in $\SSigmaZ{h}$, we need a similar result on the standard DDR space with boundary conditions
\[
\XrotZ{h}\coloneq\left\{
\uvec{v}_h\in\Xrot{h}\st\text{$v_E = 0$ for all $E\in\Ehb$}
\right\}.
\]
The interpolator $\Irot{h}:\bvec{H}^1(\Omega;\Real^2)\to\Xrot{h}$ on this space is obtained, as usual, taking $L^2$-orthogonal projections component-wise:
For all $\bvec{v}\in\bvec{H}^1(\Omega;\Real^2)$,
\[
\Irot{h}\bvec{v}
\coloneq\big(
(\Rproj{k-1}{T}\bvec{v},\cRproj{k}{T}\bvec{v})_{T\in\Th},
(\lproj{k}{E}(\bvec{v}\cdot\tangent_E))_{E\in\Eh}
\big).
\]
Using trace inequalities, one can prove the following boundedness property:
\begin{equation}\label{eq:Irot:boundedness}
  \norm[\ROT,h]{\Irot{h}\bvec{v}}\lesssim\norm[\bvec{H}^1(\Omega;\Real^2)]{\bvec{v}}
  \qquad
  \forall\bvec{v}\in\bvec{H}^1(\Omega;\Real^2).
\end{equation}

Observing that the local scalar rotor defined by \eqref{eq:RT} only depends on the unknowns in $\Xrot{T}\hookrightarrow\SSigmaZ{T}$ (see Remark \ref{rem:comparison.Xcurl}), with a little abuse of notation we use the same symbol for $\RT:\Xrot{T}\to\Poly{k}(T)$ throughout this subsection.
A similar abuse of notation is commited for $\Rh:\Xrot{h}\to\Poly{k}(\Th)$ (cf. \eqref{eq:Rh}).
The following commutation property is a consequence of \eqref{eq:RT:commutation}:
\begin{equation}\label{eq:commutation:Rh}
  \Rh\Irot{h}\bvec{v} = \lproj{k}{h}(\ROT\bvec{v})
  \qquad\forall\bvec{v}\in\bvec{H}^1(\Omega;\Real^2).
\end{equation}

\begin{lemma}[Poincar\'e inequality in $\XrotZ{h}$]\label{lem:poincare:XrotZh}
  Denote by $[\cdot,\cdot]_{\ROT,h}$ an inner product in $\XrotZ{h}$ with induced norm equivalent to $\tnorm[\ROT,h]{{\cdot}}$ uniformly in $h$.
  Then, for all $\uvec{v}_h\in\XrotZ{h}$ such that
  \begin{equation}\label{eq:poincare:XrotZh:orthogonality}
    [\uvec{v}_h, \uvec{w}_h]_{\ROT,h} = 0\qquad
    \forall\Ker\Rh\coloneq\big\{
    \uvec{w}_h\in\XrotZ{h}\st\Rh\uvec{w}_h = 0
    \big\},
  \end{equation}
  it holds,
  \[
  \tnorm[\ROT,h]{\uvec{v}_h}
  \lesssim\norm[L^2(\Omega)]{\Rh\uvec{v}_h}.
  \]
\end{lemma}

\begin{proof}
  Let $\uvec{v}_h\in\left(\Ker\Rh\right)^\bot$, where $\left(\Ker\Rh\right)^\bot$ denotes the orthogonal complement of $\Ker\Rh$ in $\XrotZ{h}$ with respect to the inner product $[\cdot,\cdot]_{\ROT,h}$.
  Owing to the surjectivity of $\ROT:\bvec{H}_0^1(\Omega;\Real^2)\to L^2_0(\Omega)$ (the proof of which is the same as for the divergence acting between the same spaces; cf., e.g., \cite{Girault.Raviart:86,Bogovskii:80,Solonnikov:01,Duran.Muschietti:01}),
  to the commutation property \eqref{eq:commutation:Rh},
  and to the uniform boundedness \eqref{eq:Irot:boundedness} of $\Irot{h}$, we infer the existence of $\bvec{v}\in\bvec{H}_0^1(\Omega;\Real^2)$ such that
  \begin{equation}\label{eq:poincare:Xrot:lifting}
  \text{
    $\Rh\uvec{v}_h = \ROT\bvec{v} = \Rh\Irot{h}\bvec{v}$
    and
    $\tnorm[\ROT,h]{\Irot{h}\bvec{v}}
    \lesssim\norm[\bvec{H}^1(\Omega;\Real^2)]{\bvec{v}}
    \lesssim\norm[L^2(\Omega)]{\Rh\uvec{v}_h}$.
  }
  \end{equation}
  We therefore have that $\uvec{v}_h - \Irot{h}\bvec{v}\in\Ker\Rh$, so that $\uvec{v}_h$ can be regarded as the $[\cdot,\cdot]_{\ROT,h}$-orthogonal projection of $\Irot{h}\bvec{v}$ on the space $\left(\Ker\Rh\right)^{\bot}$.
  Thus, by continuity of the $L^2$-orthogonal projector, the norm induced by $[\cdot,\cdot]_{\ROT,h}$ of $\uvec{v}_h$ is bounded by that of $\Irot{h}\bvec{v}$ up to a multiplicative constant independent of $h$.
  The assumed uniform equivalence between this induced norm and $\tnorm[\ROT,h]{{\cdot}}$ along with the inequality in \eqref{eq:poincare:Xrot:lifting} yields the result.
\end{proof}

\subsection{Poincar\'e inequalities in $\SSigmaZ{h}$}

\subsubsection{Poincar\'e inequality for the scalar rotor}

\begin{lemma}[Poincar\'e inequality for the scalar rotor]
  Denote by $[\cdot,\cdot]_{\ROT, h}$ an inner product for $\XrotZ{h}$ as in Lemma \ref{lem:poincare:XrotZh}.
  Then, for all $\uvec{v}_h\in\SSigmaZ{h}$ such that $\Resrot\uvec{v}_h$ satisfies \eqref{eq:poincare:XrotZh:orthogonality},
  it holds,
  \begin{equation}\label{eq:poincare:SSigmaZh}
    \tnorm[\bvec{\Sigma},h]{\uvec{v}_h}
    \lesssim\tnorm[W,h]{\uRh\uvec{v}_h}.
  \end{equation}
\end{lemma}

\begin{proof}
  By Lemma \ref{lem:poincare:XrotZh}, $\tnorm[\ROT,h]{\Resrot\uvec{v}_h}\lesssim\norm[L^2(\Omega)]{\Rh\uvec{v}_h}$.
  The result follows recalling the definition of $\tnorm[\bvec{\Sigma},h]{\uvec{v}_h}$ (see, in particular, \eqref{eq:tnorm.Sigma.T}) and noticing that the norms of the edge and vertex components that are not contained in $\tnorm[\ROT,h]{\Resrot\uvec{v}_h}$ can be bounded by the norms of the corresponding components in $\norm[W,h]{\uRh\uvec{v}_h}$ after recalling the definition \eqref{eq:uRh} of $\uRh\uvec{v}_h$ and noticing that $h_T\lesssim 1$ for all $T\in\Th$.
\end{proof}

\subsubsection{Poincar\'e inequality for rot-rot}

To state the second relevant Poincar\'e inequality in $\SSigmaZ{h}$, we need a discrete vector rotor and a scalar potential acting on the image space of $\uRh$, namely $\SWZ{h}$.
The construction presented below is similar to the one for $\SV{h}$ given in Section \ref{sec:GT.PVT}.
Specifically, we let $\VRT:\SW{T}\to\vPoly{k}(T;\Real^2)$ be such that, for all $\underline{r}_T\in\SW{T}$,
\begin{equation}\label{eq:VRT}
  \int_T\VRT\underline{r}_T\cdot\bvec{v}
  = \int_T r_T~\ROT\bvec{v}
  + \sum_{E\in\ET}\omega_{TE}\int_E \PWE\underline{r}_E~(\bvec{v}\cdot\tangent_E)
  \qquad\forall\bvec{v}\in\vPoly{k}(T;\Real^2),
\end{equation}
where, for all $E\in\ET$,
\begin{equation}\label{eq:PWE}
  \text{
    $\PWE\underline{r}_E$ is such that $\PWE\underline{r}_E(\bvec{x}_\nu) = r_\nu$ for all $\nu\in\VE$ and $\lproj{k-1}{E}\PWE\underline{r}_E = r_E$.
  }
\end{equation}
A global discrete vector rotor $\VRh:\SW{h}\to\vPoly{k}(\Th;\Real^2)$ is obtained setting, for all $\underline{r}_h\in\SW{h}$,
\[
(\VRh\underline{r}_h)_{|T} \coloneq \VRT\underline{r}_T
\qquad\forall T\in\Th.
\]

\begin{remark}[Link between the vector rotor and the discrete gradient]
  The vector rotor defined by \eqref{eq:VRT} formally coincides with the usual DDR element gradient \eqref{eq:GT} rotated by an angle of $-\nicefrac{\pi}{2}$ radians, the sole difference lying in the the degree of the element component $r_T$.
\end{remark}

We next introduce the element scalar potential $\PWT:\SW{T}\to\Poly{k+1}(T)$ such that, for all $\underline{r}_T\in\SW{T}$,
\begin{equation}\label{eq:PWT}
  \int_T\PWT\underline{r}_T~\ROT\bvec{v}
  = \int_T\VRT\underline{r}_T\cdot\bvec{v}
  - \sum_{E\in\ET}\omega_{TE}\int_E\PWE\underline{r}_E~(\bvec{v}\cdot\tangent_E)
  \qquad\forall\bvec{v}\in\cGoly{k+2}(T).
\end{equation}

\begin{remark}[Validity of \eqref{eq:PWT}]\label{rem:validity.PWT}
  Taking $\bvec{v}\in\Goly{k}(T)$ in \eqref{eq:PWT} and using the fact that $\ROT\GRAD = 0$ along with the definition \eqref{eq:VRT} of $\VRT$, it can be checked that both sides vanish, showing that \eqref{eq:PWT} holds in fact for all $\bvec{v}\in\Goly{k}(T) + \cGoly{k+2}(T)$ (hence, recalling the second decomposition in \eqref{eq:vPoly:decomposition} and the fact that the complements are hierarchical, for all $\bvec{v}\in\vPoly{k}(T;\Real^2)$).
\end{remark}

\begin{remark}[Control of the vector rotor of the element component and of the scalar potential]
  Integrating by parts the first term in the right-hand side of \eqref{eq:VRT}, we obtain
  \[
  \int_T\VROT r_T\cdot\bvec{v}
  = \int_T\VRT\underline{r}_T\cdot\bvec{v}
  + \sum_{E\in\ET}\omega_{TE}\int_E (r_T - \PWE\underline{r}_E)~(\bvec{v}\cdot\tangent_E)
  \qquad\forall\bvec{v}\in\vPoly{k}(T;\Real^2).
  \]
  Taking $\bvec{v} = \VROT r_T$, using Cauchy--Schwarz and trace inequalities, and simplifying, we infer
  \begin{equation}\label{eq:bound.VROT}
    \norm[\bvec{L}^2(T;\Real^2)]{\VROT r_T}
    \lesssim\norm[\bvec{L}^2(T;\Real^2)]{\VRT\underline{r}_T}
    + h_T^{-\nicefrac12}\sum_{E\in\ET}\norm[L^2(E)]{r_T - \PWE\underline{r}_E}.
  \end{equation}

  In a similar way, writing \eqref{eq:PWT} for $\bvec{v}\in\vPoly{k}(T;\Real^2)$ (recall Remark \ref{rem:validity.PWT}) and integrating by parts the left-hand side, we get
  \[
  \int_T\VROT\PWT\underline{r}_T\cdot\bvec{v}
  = \int_T\VRT\underline{r}_T\cdot\bvec{v}
  + \sum_{E\in\ET}\omega_{TE}\int_E(\PWT\underline{r}_T - \PWE\underline{r}_E)~(\bvec{v}\cdot\tangent_E).
  \]
  Taking $\bvec{v} = \VROT\PWT\underline{r}_T$, using Cauchy--Schwarz and trace inequalities, and simplifying, we get
  \begin{equation}\label{eq:bound.rot.PWT}
    \norm[\bvec{L}^2(T;\Real^2)]{\VROT\PWT\underline{r}_T}
    \lesssim\norm[\bvec{L}^2(T;\Real^2)]{\VRT\underline{r}_T}
    + h_T^{-\nicefrac12}\sum_{E\in\ET}\norm[L^2(E)]{\PWT\underline{r}_T - \PWE\underline{r}_E}.
  \end{equation}
\end{remark}

\begin{lemma}[Poincar\'e inequality for rot-rot]\label{lem:poincare:SSigmaZh:VROT}
  Let $\uvec{v}_h\in\SSigmaZ{h}$ be such that $\Resrot\uvec{v}_h$ satisfies \eqref{eq:poincare:XrotZh:orthogonality}.
  Then, it holds,
  \begin{equation}\label{eq:poincare:SSigmaZh:VROT}
    \tnorm[\bvec{\Sigma},h]{\uvec{v}_h}
    \lesssim\tnorm[W,h]{\uRh\uvec{v}_h}
    \lesssim\norm[\VROT\ROT,h]{\uvec{v}_h},
  \end{equation}
  where $\norm[\VROT\ROT,h]{\uvec{v}_h}^2\coloneq\sum_{T\in\Th}\norm[\VROT\ROT,T]{\uvec{v}_T}^2$ with, for all $T\in\Th$,
  \begin{multline}\label{eq:tnorm.rot.rot}
    \norm[\VROT\ROT,T]{\uvec{v}_T}^2
    \coloneq
    \norm[\bvec{L}^2(T;\Real^2)]{\VRT\uRT\uvec{v}_T}^2
    + h_T^{-2}\norm[L^2(T)]{\lproj{k}{T}\PWT\uRT\uvec{v}_T - \RT\uvec{v}_T}^2
    \\
    + h_T^{-1}\sum_{E\in\ET}\norm[L^2(E)]{\PWT\uRT\uvec{v}_T - \PWE\uRE\uvec{v}_E}^2.
  \end{multline}
\end{lemma}

\begin{proof}
  By \eqref{eq:poincare:SSigmaZh}, it holds
  \begin{equation}\label{eq:poincare:rot-rot:basic}
    \tnorm[\bvec{\Sigma},h]{\uvec{v}_h}^2
    \lesssim\tnorm[W,h]{\uRh\uvec{v}_h}^2
    \stackrel{\eqref{eq:tnorm.W.h}}{=}
    \underbrace{%
      \vphantom{\sum_{T\in\Th}\left(    
        h_T\sum_{E\in\ET}\norm[L^2(E)]{C_{\bvec{v},E}}^2
        + h_T^2\sum_{\nu\in\VT}|C_{\bvec{v},\nu}|^2.
        \right)}
      \norm[L^2(\Omega)]{\Rh\uvec{v}_h}^2
    }_{\term_1}
    + \underbrace{%
      \sum_{T\in\Th}\left(    
      h_T\sum_{E\in\ET}\norm[L^2(E)]{C_{\bvec{v},E}}^2
      + h_T^2\sum_{\nu\in\VT}|C_{\bvec{v},\nu}|^2
      \right)
    }_{\term_2}.x
  \end{equation}
  We next proceed to bound $\term_1$ and $\term_2$ by $\norm[\VROT\ROT,h]{\uvec{v}_h}$.
  \medskip\\
  \underline{(i) \emph{Estimate of $\term_1$.}}
  Applying the discrete Poincaré inequality in HHO spaces \cite[Lemma 2.15]{Di-Pietro.Droniou:20} to the vector $\big( (\RT\uvec{v}_T)_{T\in\Th}, (\PWE\uRE\uvec{v}_E)_{E\in\Eh}\big)$  (which has vanishing boundary components since $\uvec{v}_h\in\SSigmaZ{h}$, as can be checked writing \eqref{eq:PWE} with $\underline{r}_E = \uRE\uvec{v}_E$), we get
  \[
  \norm[L^2(\Omega)]{\Rh\uvec{v}_h}^2
  \lesssim
  \sum_{T\in\Th}\left(
  \norm[\bvec{L}^2(T;\Real^2)]{\VROT\RT\uvec{v}_T}^2
  + h_T^{-1}\sum_{E\in\ET}\norm[L^2(E)]{\RT\uvec{v}_T - \PWE\uRE\uvec{v}_E}^2
  \right),
  \]
  where we have additionally used the fact that
  \begin{equation}\label{eq:norm.grad.rot}
    \norm[\bvec{L}^2(T;\Real^2)]{\GRAD\cdot}
    = \norm[\bvec{L}^2(T;\Real^2)]{\rotation{-\nicefrac{\pi}{2}}\GRAD\cdot}
    = \norm[\bvec{L}^2(T;\Real^2)]{\VROT\cdot}
  \end{equation}
  for the first term in the right-hand side.
  The estimate on $\term_1$ follows provided we show that
  \begin{equation}\label{eq:poincare:rot-rot:estimate.rhs}
    \sum_{T\in\Th}\left(
    \norm[\bvec{L}^2(T;\Real^2)]{\VROT\RT\uvec{v}_T}^2
    + h_T^{-1}\sum_{E\in\ET}\norm[L^2(E)]{\RT\uvec{v}_T - \PWE\uRE\uvec{v}_E}^2
    \right)
    \lesssim\norm[\VROT\ROT,h]{\uvec{v}_h}^2.
  \end{equation}
  Writing \eqref{eq:bound.VROT} with $\underline{r}_T = \uRT\uvec{v}_T$, squaring, and using the fact that $\card(\ET)\lesssim 1$ by mesh regularity, we get
  \begin{multline} \label{eq:poincare:rot-rot:basic}
    \sum_{T\in\Th}\left(
    \norm[\bvec{L}^2(T;\Real^2)]{\VROT\RT\uvec{v}_T}^2
    + h_T^{-1}\sum_{E\in\ET}\norm[L^2(E)]{\RT\uvec{v}_T - \PWE\uRE\uvec{v}_E}^2
    \right)
    \\
    \lesssim
    \sum_{T\in\Th}\left(
    \norm[\bvec{L}^2(T;\Real^2)]{\VRT\RT\uvec{v}_T}
    + h_T^{-1}\sum_{E\in\ET}\norm[L^2(E)]{\RT\uvec{v}_T - \PWE\uRE\uvec{v}_E}^2
    \right).
  \end{multline}
  It only remains to estimate the boundary term.
  To this end, for a given $E\in\ET$, we insert $\pm(\PWT\uRT\uvec{v}_T - \lproj{k}{T}\PWT\uRT\uvec{v}_T)$ into the norm and use triangle and trace inequalities to write
  \[
  \begin{aligned}
    \norm[L^2(E)]{\RT\uvec{v}_T - \PWE\uRE\uvec{v}_E}
    &\le h_T^{-\nicefrac12}\norm[L^2(T)]{\lproj{k}{T}\PWT\uRT\uvec{v}_T - \RT\uvec{v}_T}
    \\
    &\quad
    + h_T^{-\nicefrac12}\norm[L^2(T)]{\PWT\uRT\uvec{v}_T - \lproj{k}{T}\PWT\uRT\uvec{v}_T}    
    \\
    &\quad
    + \norm[L^2(E)]{\PWT\uRT\uvec{v}_T - \PWE\uRE\uvec{v}_E}.
  \end{aligned}
  \]
  Using the approximation properties of $\lproj{k}{T}$ (cf. \cite{Di-Pietro.Droniou:17} and \cite[Chapter~1]{Di-Pietro.Droniou:20} concerning the extension to non-star-shaped elements), we next write
  \[
  \begin{aligned}
    &\norm[L^2(T)]{\PWT\uRT\uvec{v}_T - \lproj{k}{T}\PWT\uRT\uvec{v}_T}
    \\
    &\quad\lesssim
    h_T\norm[\bvec{L}^2(T;\Real^2)]{\GRAD\PWT\uRT\uvec{v}_T}
    \stackrel{\eqref{eq:norm.grad.rot}}{=}
    h_T\norm[\bvec{L}^2(T;\Real^2)]{\VROT\PWT\uRT\uvec{v}_T}
    \\
    &\quad
    \!\stackrel{\eqref{eq:bound.rot.PWT}}{\lesssim}
    h_T\norm[\bvec{L}^2(T;\Real^2)]{\VRT\uRT\uvec{v}_T}
    + h_T^{\nicefrac12}\sum_{F\in\ET}\norm[L^2(F)]{\PWT\uRT\uvec{v}_T - \PWF\uRF\uvec{v}_F},
  \end{aligned}
  \]
  so that
  \begin{equation}\label{eq:poincare:rot-rot:1}
    \begin{aligned}
      \norm[L^2(E)]{\RT\uvec{v}_T - \PWE\uRE\uvec{v}_E}
      &\lesssim
      h_T^{-\nicefrac12}\norm[L^2(T)]{\lproj{k}{T}\PWT\uRT\uvec{v}_T - \RT\uvec{v}_T}
      \\
      &\quad
      + h_T^{\nicefrac12}\norm[\bvec{L}^2(T;\Real^2)]{\VRT\uRT\uvec{v}_T}
      \\
      &\quad
      + \sum_{F\in\ET}\norm[L^2(F)]{\PWT\uRT\uvec{v}_T - \PWF\uRF\uvec{v}_F}.
  \end{aligned}
  \end{equation}
  Plugging \eqref{eq:poincare:rot-rot:1} into \eqref{eq:poincare:rot-rot:basic} and using the fact that $\card(\ET)\lesssim 1$, it follows that
  \begin{equation}\label{eq:poincare:rot-rot:T1}
    \term_1\lesssim\norm[\VROT\ROT,h]{\uvec{v}_h}.
  \end{equation}
  \noindent\underline{(ii) \emph{Estimate of $\term_2$.}}
  We start by recalling the following norm equivalence, valid for all $E\in\Eh$, all polynomial degrees $m\ge 0$, and all $\varphi\in\Poly{m}(E)$:
  \begin{equation}\label{eq:norm.equivalence:1d}
    \norm[L^2(E)]{\varphi}^2
    \simeq\norm[L^2(E)]{\lproj{m-2}{E}\varphi}^2
    + h_E\sum_{\nu\in\VE}\varphi(\bvec{x}_\nu)^2
  \end{equation}
  The proof of this equivalence follows recalling the estimate of the $L^2$-norm of functions on the unit segment $[0,1]$ corresponding to the first display equation in the proof in Eq. (58) of \cite[Proposition 10]{Di-Pietro.Droniou:23} and using the isomorphism $[0,1]\ni s\mapsto\bvec{x}_{V_1} + s h_E(\bvec{x} - \bvec{x}_{V_1})$ (with $V_1$ denoting the first vertex of $E$ in the direction of $\tangent_E$).
  Using \eqref{eq:norm.equivalence:1d} with $(\varphi,m) = (\PWE\uRE\uvec{v}_E,k+1)$ and recalling that, by \eqref{eq:PWE} written for $\underline{r}_E = \uRE\uvec{v}_E$, $\lproj{k-1}{E}\PWE\uRE\uvec{v}_E = C_{\bvec{v},E}$ and $\PWE\uRE\uvec{v}_E(\bvec{x}_\nu) = C_{\bvec{v},\nu}$ for all $\nu\in\VE$, we get
  \[
  \begin{aligned}
    \term_2
    &\lesssim \sum_{T\in\Th}h_T\sum_{E\in\ET}\norm[L^2(E)]{\PWE\uRE\uvec{v}_E}^2
    \\
    &\lesssim\sum_{T\in\Th}\norm[L^2(T)]{\PWT\uRT\uvec{v}_T}^2
    + \sum_{T\in\Th}h_T\sum_{E\in\ET}\norm[L^2(E)]{\PWT\uRT\uvec{v}_T - \PWE\uRE\uvec{v}_E}^2
    \eqcolon \term_{2,1} + \term_{2,2},
  \end{aligned}
  \]
  where, to pass to the second line, we have inserted $\pm\PWT\uRT\uvec{v}_T$ into the norm and used triangle and discrete trace inequalities.
  Applying again \cite[Lemma 2.15]{Di-Pietro.Droniou:20}, this time to $\big( (\PWT\uRT\uvec{v}_T)_{T\in\Th}, (\PWE\uRE\uvec{v}_E)_{E\in\Eh}\big)$, and using \eqref{eq:norm.grad.rot}, we get
  \[
  \term_{2,1}
  \lesssim\sum_{T\in\Th}\left(
  \norm[\bvec{L}^2(T;\Real^2)]{\VROT\PWT\uRT\uvec{v}_T}^2
  + h_T^{-1}\sum_{E\in\ET}\norm[L^2(E)]{\PWT\uRT\uvec{v}_T - \PWE\uRE\uvec{v}_E}^2
  \right),
  \]
  hence $\term_{2,1}\lesssim\norm[\VROT\ROT,h]{\uvec{v}_h}^2$ by \eqref{eq:bound.rot.PWT} with $\underline{r}_T = \uRh\uvec{v}_T$ for all $T\in\Th$.
  On the other hand, using $h_T\lesssim 1$ along with the definition \eqref{eq:tnorm.rot.rot} of $\norm[\VROT\ROT,h]{{\cdot}}$, we also have $\term_{2,2}\lesssim\norm[\VROT\ROT,h]{\uvec{v}_h}^2$.
  Combining these estimates gives
  \begin{equation}\label{eq:poincare:rot-rot:T2}
    \term_2\lesssim\norm[\VROT\ROT,h]{\uvec{v}_h}^2.
  \end{equation}
  \noindent\underline{(iii) \emph{Conclusion}.}
  Plug \eqref{eq:poincare:rot-rot:T1} and \eqref{eq:poincare:rot-rot:T2} into \eqref{eq:poincare:rot-rot:basic} and take the square root of the resulting inequality.
\end{proof}

%------------------------------------------------------------------------------%

\section{Application to a quad-rot problem}\label{sec:application}

Denote by $\normal$ the unit outward normal vector field on the boundary $\partial\Omega$ of $\Omega$, and set $\HrotrotZ{\Omega}\coloneq\left\{\bvec{v}\in\Hrotrot{\Omega}\st\text{$\normal\times\bvec{v} = \bvec{0}$ and $\ROT\bvec{u} = 0$ on $\partial\Omega$}\right\}$.
As customary, we also denote by $H_0^1(\Omega)$ the subspace of $H^1(\Omega)$ spanned by functions with vanishing scalar trace on $\partial\Omega$.
We use the discrete complex developed in the previous sections to discretise the following problem:
Given $\bvec{f}\in\bvec{L}^2(\Omega)$, find $\bvec{u}\in\HrotrotZ{\Omega}$ and $p\in H^1_0(\Omega)$ such that
\begin{equation}\label{eq:weak}
  \begin{alignedat}{4}
    a(\bvec{u}, \bvec{v}) + b(\bvec{v}, p) &= \int_\Omega\bvec{f}\cdot\bvec{v}
    &\qquad&\forall\bvec{v}\in\HrotrotZ{\Omega},
    \\
    -b(\bvec{u}, q) &= 0
    &\qquad&\forall q\in H_0^1(\Omega),
  \end{alignedat}
\end{equation}
with bilinear forms $a:\left[\HrotrotZ{\Omega}\right]^2\to\Real$ and $b:\HrotrotZ{\Omega}\times H_0^1(\Omega)\to\Real$ such that, for all $(\bvec{w}, \bvec{v}, q)\in\HrotrotZ{\Omega}\times\HrotrotZ{\Omega}\times H_0^1(\Omega)$,
\[
a(\bvec{w}, \bvec{v})\coloneq\int_\Omega\VROT\ROT\bvec{w}\cdot\VROT\ROT\bvec{v},\qquad
b(\bvec{v}, q)\coloneq\int_\Omega\bvec{v}\cdot\GRAD q.
\]
The solution to this problem satisfies almost everywhere
\[
\begin{alignedat}{4}
  (\VROT\ROT)^2\bvec{u} + \GRAD p &= \bvec{f}
  &\qquad&\text{in $\Omega$},
  \\
  \DIV\bvec{u} &= 0
  &\qquad&\text{in $\Omega$},
  \\  
  \bvec{u}\times\normal &= 0
  &\qquad&\text{on $\partial\Omega$},
  \\
  \ROT\bvec{u} &= 0
  &\qquad&\text{on $\partial\Omega$}.
\end{alignedat}
\]

\begin{remark}[Forcing term]\label{eq:forcing.term}
  Notice that we consider here a more realistic formulation with respect to other works in the literature, since the forcing term is not prepared, i.e., it is not divergence-free nor expressed as $\VROT\varphi$ in terms of a scalar potential $\varphi$.
  This implies, in particular, that the Lagrange multiplier $p$ can be non-zero.
  Moreover, the fact that we do not assume $\varphi$ available prevents us from discretising the right-hand based on the form $\int_\Omega\varphi~\ROT\bvec{v}$ (as done, e.g., in \cite{Zhao.Zhang:21}) instead of $\int_F\bvec{f}\cdot\bvec{v}$.
  As a result, an adjoint consistency property for rot-rot is required, which makes the object of Lemma \ref{lem:adjoint.consistency:rot-rot} below.
  This result appears to be the first of this kind for polygonal methods.
\end{remark}

\subsection{Discrete problem and main results}

We define the global discrete bilinear forms $a_h:\SSigma{h}\times\SSigma{h}\to\Real$, $b_h:\SSigma{h}\times\SV{h}\to\Real$,  and $\ell_h:\bvec{L}^2(\Omega;\Real^2)\times\SSigma{h}\to\Real$ such that, for all $(\uvec{w}_h,\uvec{v}_h,\underline{q}_h,\bvec{g})\in\SSigma{h}\times\SSigma{h}\times\SV{h}\times\bvec{L}^2(\Omega;\Real^2)$,
\begin{equation}\label{eq:ah.bh.lh}
  \begin{gathered}
    a_h(\uvec{w}_h,\uvec{v}_h)
    \coloneq\sum_{T\in\Th}a_T(\uvec{w}_T,\uvec{v}_T),\qquad
    b_h(\uvec{v}_h,\underline{q}_h)
    \coloneq(\uvec{v}_h,\uGh\underline{q}_h)_{\btens{\Sigma},h},
    \\
    \ell_h(\bvec{g}, \uvec{v}_h)\coloneq\sum_{T\in\Th}\int_T\bvec{f}\cdot\PSigmaT\uvec{v}_T,
  \end{gathered}
\end{equation}
where, for all $T\in\Th$, the local discrete bilinear form $a_T:\SSigma{T}\times\SSigma{T}\to\Real$ is such that
\begin{equation}\label{eq:aT}
  a_T(\uvec{w}_T,\uvec{v}_T)
  \coloneq
  \int_T\VRT\uRT\uvec{w}_T\cdot\VRT\uRT\uvec{v}_T
  + s_T(\uRT\uvec{w}_T,\uRT\uvec{v}_T),
\end{equation}
with stabilisation bilinear form $s_T:\SW{T}\times\SW{T}\to\Real$ such that, for all $(\underline{r}_T, \underline{q}_T)\in\SW{T}\times\SW{T}$,
\begin{equation}\label{eq:sT}
  \begin{aligned}
    s_T(\underline{r}_T, \underline{q}_T)
    &\coloneq
    h_T^{-2}\int_T(\lproj{k}{T}\PWT\underline{r}_T - r_T)~(\lproj{k}{T}\PWT\underline{q}_T - q_T)
    \\
    &\quad
    + h_T^{-1}\sum_{E\in\ET}\int_E(\PWT\underline{r}_T - \PWE\underline{r}_E)~(\PWT\underline{q}_T - \PWE\underline{q}_E).
  \end{aligned}
\end{equation}
The discrete problem reads:
Find $(\uvec{u}_h,\underline{p}_h)\in\SSigmaZ{h}\times\SVZ{h}$ such that
\begin{equation}\label{eq:discrete}
  \begin{alignedat}{4}
    a_h(\uvec{u}_h,\uvec{v}_h) + b_h(\uvec{v}_h,\underline{p}_h)
    &= \ell_h(\bvec{f}, \uvec{v}_h)
    &\qquad&
    \forall\uvec{v}_h\in\SSigmaZ{h},
    \\
    -b_h(\uvec{u}_h,\underline{q}_h) &= 0
    &\qquad&
    \forall\underline{q}_h\in\SVZ{h}.
  \end{alignedat}
\end{equation}
The variational formulation corresponding to \eqref{eq:discrete} is:
Find $(\uvec{u}_h,\underline{p}_h)\in\SSigmaZ{h}\times\SVZ{h}$ such that
\begin{equation}\label{eq:discrete:variational}
  \mathcal{A}_h((\uvec{u}_h, \underline{p}_h), (\uvec{v}_h, \underline{q}_h))
  = \ell_h(\bvec{f}, \uvec{v}_h)
  \qquad\forall(\uvec{v}_h, \underline{q}_h)\in\SSigmaZ{h}\times\SVZ{h},
\end{equation}
with bilinear form $\mathcal{A}_h:\big[\SSigma{h}\times\SV{h}\big]^2\to\Real$ such that, for all $(\uvec{w}_h,\underline{r}_h)$ and all $(\uvec{v}_h,\underline{q}_h)$ in $\SSigma{h}\times\SV{h}$,
\begin{equation}\label{eq:Ah}
  \mathcal{A}_h((\uvec{w}_h, \underline{r}_h), (\uvec{v}_h, \underline{q}_h))
  \coloneq a_h(\uvec{w}_h,\uvec{v}_h)
  + b_h(\uvec{v}_h,\underline{r}_h)
  - b_h(\uvec{w}_h,\underline{q}_h).
\end{equation}

\begin{remark}[Serendipity version of the scheme]
  A serendipity version  of the above scheme is obtained in the spirit of \cite[Section 6.6]{Di-Pietro.Droniou:22} substituting the spaces $\SSigmaZ{h}$ and $\SVZ{h}$ with the serendipity versions defined in Section \ref{sec:discrete.complex:serendipity}.
  The analysis remains unchanged.
\end{remark}

We equip $\SSigmaZ{h}\times\SVZ{h}$ with the following graph norm:
For all $(\uvec{v}_h, \underline{q}_h)\in\SSigmaZ{h}\times\SVZ{h}$,
\begin{equation}\label{eq:norm.Sigma.x.V.h}
  \norm[\bvec{\Sigma}\times V,h]{(\uvec{v}_h, \underline{q}_h)}^2
  \coloneq
  \tnorm[\bvec{\Sigma},h]{\uvec{v}_h}^2
  + \norm[\VROT\ROT,h]{\uvec{v}_h}^2
  + \tnorm[V,h]{\underline{q}_h}^2
  + \tnorm[\bvec{\Sigma},h]{\uGh\underline{q}_h}^2.
\end{equation}
The main results for the numerical scheme \eqref{eq:discrete} (or, equivalently, \eqref{eq:discrete:variational}) are stated below.

\begin{lemma}[Inf-sup stability]\label{lem:inf-sup}
  It holds, for all $(\uvec{w}_h,\underline{r}_h)\in\SSigmaZ{h}\times\SVZ{h}$,
  \begin{equation}\label{eq:inf-sup}
    \norm[\bvec{\Sigma}\times V,h]{(\uvec{w}_h, \underline{r}_h)}
    \lesssim
    \sup_{(\uvec{v}_h,\underline{q}_h)\in\SSigmaZ{h}\times\SVZ{h}\setminus\{0\}}
    \frac{%
      \mathcal{A}_h((\uvec{w}_h, \underline{r}_h), (\uvec{v}_h, \underline{q}_h))
    }{%
      \norm[\bvec{\Sigma}\times V,h]{(\uvec{v}_h, \underline{q}_h)}
    }.
  \end{equation}
\end{lemma}

\begin{proof}
  See Section \ref{sec:application:stability}.
\end{proof}

In what follows, we denote by $H^s(\Th)$ and $\bvec{H}^s(\Th;\Real^2)$ the scalar-valued and vector-valued broken Hilbert spaces on the mesh.

\begin{theorem}[Error estimate]\label{thm:error.estimate}
  Denote by $(\bvec{u},p)\in\HrotrotZ{\Omega}\times H^1_0(\Omega)$ the unique solution to \eqref{eq:weak} and define the seminorm $|{\cdot}|_{\bvec{H}^{(k+1,2)}(\Th;\Real^2)}$ such that, for all $\bvec{v}\in\bvec{H}^{\max(k+1,2)}(\Th;\Real^2)$,
  \[
  |\bvec{v}|_{\bvec{H}^{(k+1,2)}(\Th;\Real^2)}^2\begin{cases}
  \sum_{T\in\Th}\big(
  \seminorm[\bvec{H}^1(T;\Real^2)]{\bvec{v}}^2
  + h_T^2\seminorm[\bvec{H}^2(T;\Real^2)]{\bvec{v}}^2
  \big)
  & \text{if $k=0$},
  \\
  \seminorm[\bvec{H}^{k+1}(\Th;\Real^2)]{\bvec{v}}^2
  & \text{if $k\ge 1$}.
  \end{cases}
  \]%
  Recalling the second complex of \eqref{eq:commutation}, assume the following additional regularity:
  $\bvec{u}\in\bvec{\Sigma}\cap\bvec{H}^{\max(k+1,2)}(\Th;\Real^2)$,
  $\ROT\bvec{u}\in H^{k+2}(\Th)$,
  $\VROT\ROT\bvec{u}\in\bvec{H}^1(\Omega;\Real^2)\cap\bvec{H}^{k+1}(\Th;\Real^2)$,
  $\ROT(\VROT\ROT\bvec{u})\in H^1(\Omega)\cap H^{k+2}(\Th)$,
  %% $(\VROT\ROT)^2\bvec{u}\in\bvec{L}^2(\Omega;\Real^2)$,
  and $p\in H^{k+2}(\Th)$.
  Then, denoting by $(\uvec{u}_h,\underline{p}_h)\in\SSigmaZ{h}\times\SVZ{h}$ the unique solution to \eqref{eq:discrete} (or, equivalently, \eqref{eq:discrete:variational}) it holds
  \begin{equation}\label{eq:err.estimate}
    \norm[\bvec{\Sigma}\times V,h]{(\uvec{u}_h - \ISigma{h}\bvec{u}, \underline{p}_h - \IV{h} p)}
    \lesssim h^{k+1}\left(
    \mathcal{N}(\bvec{u})
    + \seminorm[\bvec{H}^{(k+1,2)}(\Th;\Real^2)]{\bvec{u}}
    + \seminorm[H^{k+2}(\Th)]{p}
    \right),
  \end{equation}
  with
  \begin{equation}\label{eq:mathcal.N}
    \mathcal{N}(\bvec{u})\coloneq
    \seminorm[H^{k+2}(\Th)]{\ROT\bvec{u}}
    + \seminorm[\bvec{H}^{k+1}(\Th;\Real^2)]{\VROT\ROT\bvec{u}}
    + \seminorm[H^{k+2}(\Th)]{\ROT(\VROT\ROT\bvec{u})}.
  \end{equation}
\end{theorem}

\begin{proof}
  See Section \ref{sec:application:err.estimate}.
\end{proof}

\begin{remark}[Error estimate in double-bar norms]\label{rem:err.estimate:operator.norms}
  Starting from \eqref{eq:err.estimate} and using the norm equivalence \eqref{eq:norm.equivalence:SSigma.h} for $\tnorm[\bvec{\Sigma},h]{{\cdot}}$ and those resulting from \cite[Proposition 6]{Di-Pietro.Droniou:21*1} for $\tnorm[V,h]{{\cdot}}$, one can easily infer the following error estimate where (triple-bar) component norms are replaced by the corresponding (double-bar) operator norms:
  \begin{multline}\label{eq:err.estimate:operator.norms}
  \norm[\bvec{\Sigma},h]{\uvec{u}_h - \ISigma{h}\bvec{u}}
  + \norm[\VROT\ROT,h]{\uvec{u}_h - \ISigma{h}\bvec{u}}
  + \norm[V,h]{\underline{p}_h - \IV{h}p}
  + \norm[\bvec{\Sigma},h]{\uGh(\underline{p}_h - \IV{h}p)}
  \\
  \lesssim h^{k+1}\left(
    \mathcal{N}(\bvec{u})
    + \seminorm[\bvec{H}^{(k+1,2)}(\Th;\Real^2)]{\bvec{u}}
    + \seminorm[H^{k+2}(\Th)]{p}
    \right),
  \end{multline}
  where $\norm[\bvec{\Sigma},h]{{\cdot}}$ is defined by \eqref{eq:norm.Sigma}, while $\norm[V,h]{{\cdot}}$ is the standard DDR component norm for the $H^1$-like space defined in \cite[Section 4.4]{Di-Pietro.Droniou:21*1}.
\end{remark}

\subsection{Numerical results}\label{sec:application:numerical.results}

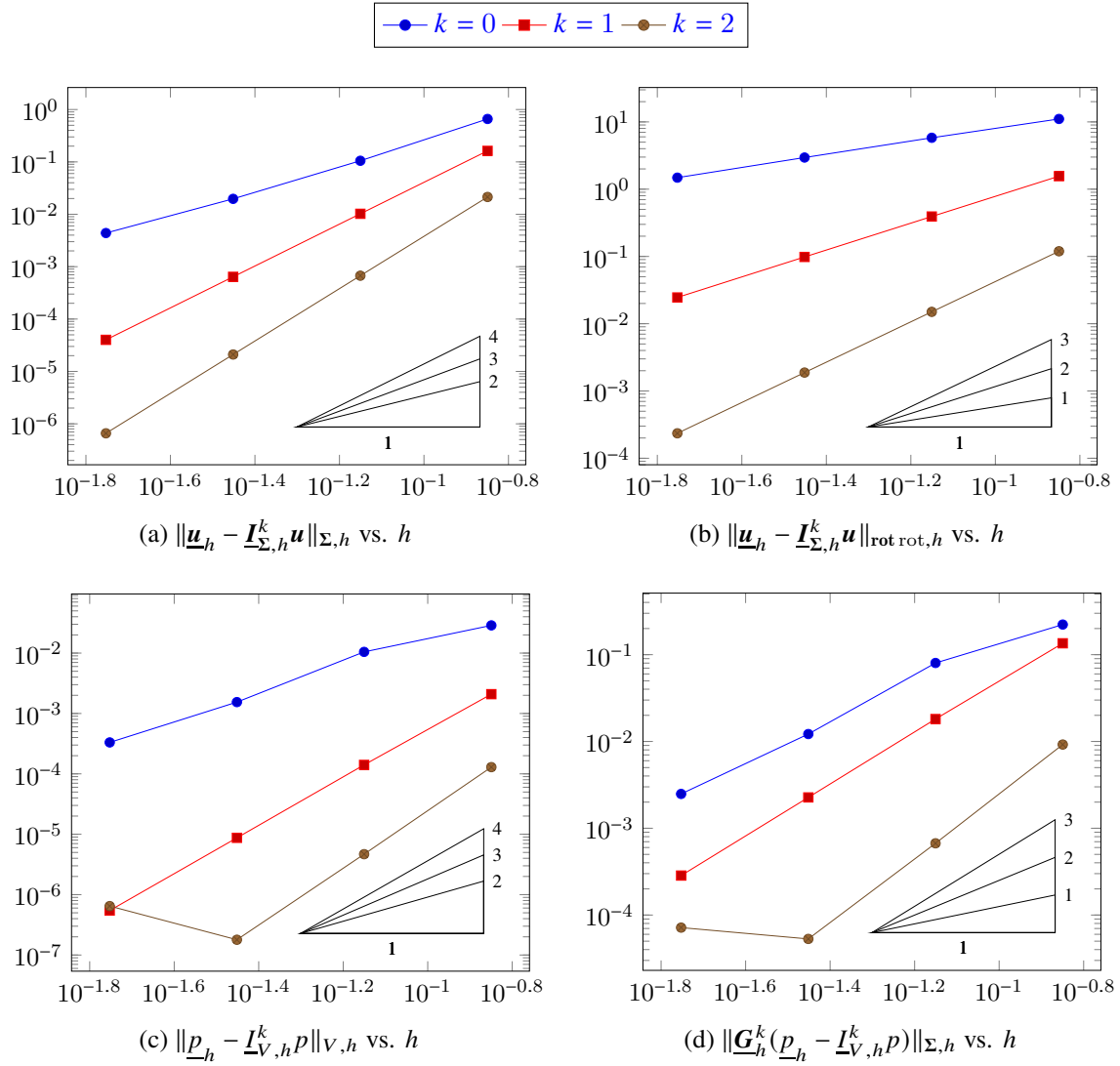
\begin{figure}\centering
  \ref{leg:cartesian}
  \vspace{0.5cm}\\
  \begin{minipage}[b]{0.475\linewidth}\centering
    \begin{tikzpicture}[scale=0.90]
      \begin{loglogaxis}[
          legend columns=-1,
          legend to name=leg:cartesian
        ]
        \addplot table[x=MeshSize,y=ErrUL2] {dat/cart_k0/data_rates.dat};
        \addplot table[x=MeshSize,y=ErrUL2] {dat/cart_k1/data_rates.dat};
        \addplot table[x=MeshSize,y=ErrUL2] {dat/cart_k2/data_rates.dat};

        \logLogSlopeTriangle{0.9}{0.4}{0.1}{2}{black};
        \logLogSlopeTriangle{0.9}{0.4}{0.1}{3}{black};
        \logLogSlopeTriangle{0.9}{0.4}{0.1}{4}{black};
        
        \legend{$k=0$,$k=1$,$k=2$};
      \end{loglogaxis}
    \end{tikzpicture}
    \subcaption{$\norm[\bvec{\Sigma},h]{\uvec{u}_h - \ISigma{h}\bvec{u}}$ vs. $h$}
  \end{minipage}
  \begin{minipage}[b]{0.475\linewidth}\centering
    \begin{tikzpicture}[scale=0.90]
      \begin{loglogaxis}
        \addplot table[x=MeshSize,y=ErrURotRot] {dat/cart_k0/data_rates.dat};
        \addplot table[x=MeshSize,y=ErrURotRot] {dat/cart_k1/data_rates.dat};
        \addplot table[x=MeshSize,y=ErrURotRot] {dat/cart_k2/data_rates.dat};

        \logLogSlopeTriangle{0.9}{0.4}{0.1}{1}{black};
        \logLogSlopeTriangle{0.9}{0.4}{0.1}{2}{black};
        \logLogSlopeTriangle{0.9}{0.4}{0.1}{3}{black};
      \end{loglogaxis}
    \end{tikzpicture}
    \subcaption{$\norm[\VROT\ROT,h]{\uvec{u}_h - \ISigma{h}\bvec{u}}$ vs. $h$}
  \end{minipage}
  \vspace{0.5cm}\\
  \begin{minipage}[b]{0.475\linewidth}\centering
    \begin{tikzpicture}[scale=0.90]
      \begin{loglogaxis}
        \addplot table[x=MeshSize,y=ErrPL2] {dat/cart_k0/data_rates.dat};
        \addplot table[x=MeshSize,y=ErrPL2] {dat/cart_k1/data_rates.dat};
        \addplot table[x=MeshSize,y=ErrPL2] {dat/cart_k2/data_rates.dat};

        \logLogSlopeTriangle{0.9}{0.4}{0.1}{2}{black};
        \logLogSlopeTriangle{0.9}{0.4}{0.1}{3}{black};
        \logLogSlopeTriangle{0.9}{0.4}{0.1}{4}{black};
      \end{loglogaxis}
    \end{tikzpicture}
    \subcaption{$\norm[V,h]{\underline{p}_h - \IV{h} p}$ vs. $h$}
  \end{minipage}
  \begin{minipage}[b]{0.475\linewidth}\centering
    \begin{tikzpicture}[scale=0.90]
      \begin{loglogaxis}
        \addplot table[x=MeshSize,y=ErrPGrad] {dat/cart_k0/data_rates.dat};
        \addplot table[x=MeshSize,y=ErrPGrad] {dat/cart_k1/data_rates.dat};
        \addplot table[x=MeshSize,y=ErrPGrad] {dat/cart_k2/data_rates.dat};

        \logLogSlopeTriangle{0.9}{0.4}{0.1}{1}{black};
        \logLogSlopeTriangle{0.9}{0.4}{0.1}{2}{black};
        \logLogSlopeTriangle{0.9}{0.4}{0.1}{3}{black};
      \end{loglogaxis}
    \end{tikzpicture}
        \subcaption{$\norm[\bvec{\Sigma},h]{\uGh(\underline{p}_h - \IV{h} p)}$ vs. $h$}
  \end{minipage}
  \caption{Convergence of the error components for the numerical test of Section \ref{sec:application:numerical.results} on the Cartesian orthogonal mesh family.\label{fig:cartesian}}
\end{figure}

\begin{figure}\centering
  \ref{leg:triangular}
  \vspace{0.5cm}\\
  \begin{minipage}[b]{0.475\linewidth}\centering
    \begin{tikzpicture}[scale=0.90]
      \begin{loglogaxis}[
          legend columns=-1,
          legend to name=leg:triangular
        ]
        \addplot table[x=MeshSize,y=ErrUL2] {dat/tri_k0/data_rates.dat};
        \addplot table[x=MeshSize,y=ErrUL2] {dat/tri_k1/data_rates.dat};
        \addplot table[x=MeshSize,y=ErrUL2] {dat/tri_k2/data_rates.dat};

        \logLogSlopeTriangle{0.9}{0.4}{0.1}{2}{black};
        \logLogSlopeTriangle{0.9}{0.4}{0.1}{3}{black};
        \logLogSlopeTriangle{0.9}{0.4}{0.1}{4}{black};
        
        \legend{$k=0$,$k=1$,$k=2$};
      \end{loglogaxis}
    \end{tikzpicture}
    \subcaption{$\norm[\bvec{\Sigma},h]{\uvec{u}_h - \ISigma{h}\bvec{u}}$ vs. $h$}
  \end{minipage}
  \begin{minipage}[b]{0.475\linewidth}\centering
    \begin{tikzpicture}[scale=0.90]
      \begin{loglogaxis}
        \addplot table[x=MeshSize,y=ErrURotRot] {dat/tri_k0/data_rates.dat};
        \addplot table[x=MeshSize,y=ErrURotRot] {dat/tri_k1/data_rates.dat};
        \addplot table[x=MeshSize,y=ErrURotRot] {dat/tri_k2/data_rates.dat};

        \logLogSlopeTriangle{0.9}{0.4}{0.1}{1}{black};
        \logLogSlopeTriangle{0.9}{0.4}{0.1}{2}{black};
        \logLogSlopeTriangle{0.9}{0.4}{0.1}{3}{black};
      \end{loglogaxis}
    \end{tikzpicture}
    \subcaption{$\norm[\VROT\ROT,h]{\uvec{u}_h - \ISigma{h}\bvec{u}}$ vs. $h$}
  \end{minipage}
  \vspace{0.5cm}\\
  \begin{minipage}[b]{0.475\linewidth}\centering
    \begin{tikzpicture}[scale=0.90]
      \begin{loglogaxis}
        \addplot table[x=MeshSize,y=ErrPL2] {dat/tri_k0/data_rates.dat};
        \addplot table[x=MeshSize,y=ErrPL2] {dat/tri_k1/data_rates.dat};
        \addplot table[x=MeshSize,y=ErrPL2] {dat/tri_k2/data_rates.dat};

        \logLogSlopeTriangle{0.9}{0.4}{0.1}{2}{black};
        \logLogSlopeTriangle{0.9}{0.4}{0.1}{3}{black};
        \logLogSlopeTriangle{0.9}{0.4}{0.1}{4}{black};
      \end{loglogaxis}
    \end{tikzpicture}
    \subcaption{$\norm[V,h]{\underline{p}_h - \IV{h} p}$ vs. $h$}
  \end{minipage}
  \begin{minipage}[b]{0.475\linewidth}\centering
    \begin{tikzpicture}[scale=0.90]
      \begin{loglogaxis}
        \addplot table[x=MeshSize,y=ErrPGrad] {dat/tri_k0/data_rates.dat};
        \addplot table[x=MeshSize,y=ErrPGrad] {dat/tri_k1/data_rates.dat};
        \addplot table[x=MeshSize,y=ErrPGrad] {dat/tri_k2/data_rates.dat};

        \logLogSlopeTriangle{0.9}{0.4}{0.1}{1}{black};
        \logLogSlopeTriangle{0.9}{0.4}{0.1}{2}{black};
        \logLogSlopeTriangle{0.9}{0.4}{0.1}{3}{black};
      \end{loglogaxis}
    \end{tikzpicture}
        \subcaption{$\norm[\bvec{\Sigma},h]{\uGh(\underline{p}_h - \IV{h} p)}$ vs. $h$}
  \end{minipage}
  \caption{Convergence of the error components for the numerical test of Section \ref{sec:application:numerical.results} on the triangular mesh family.\label{fig:triangular}}
\end{figure}
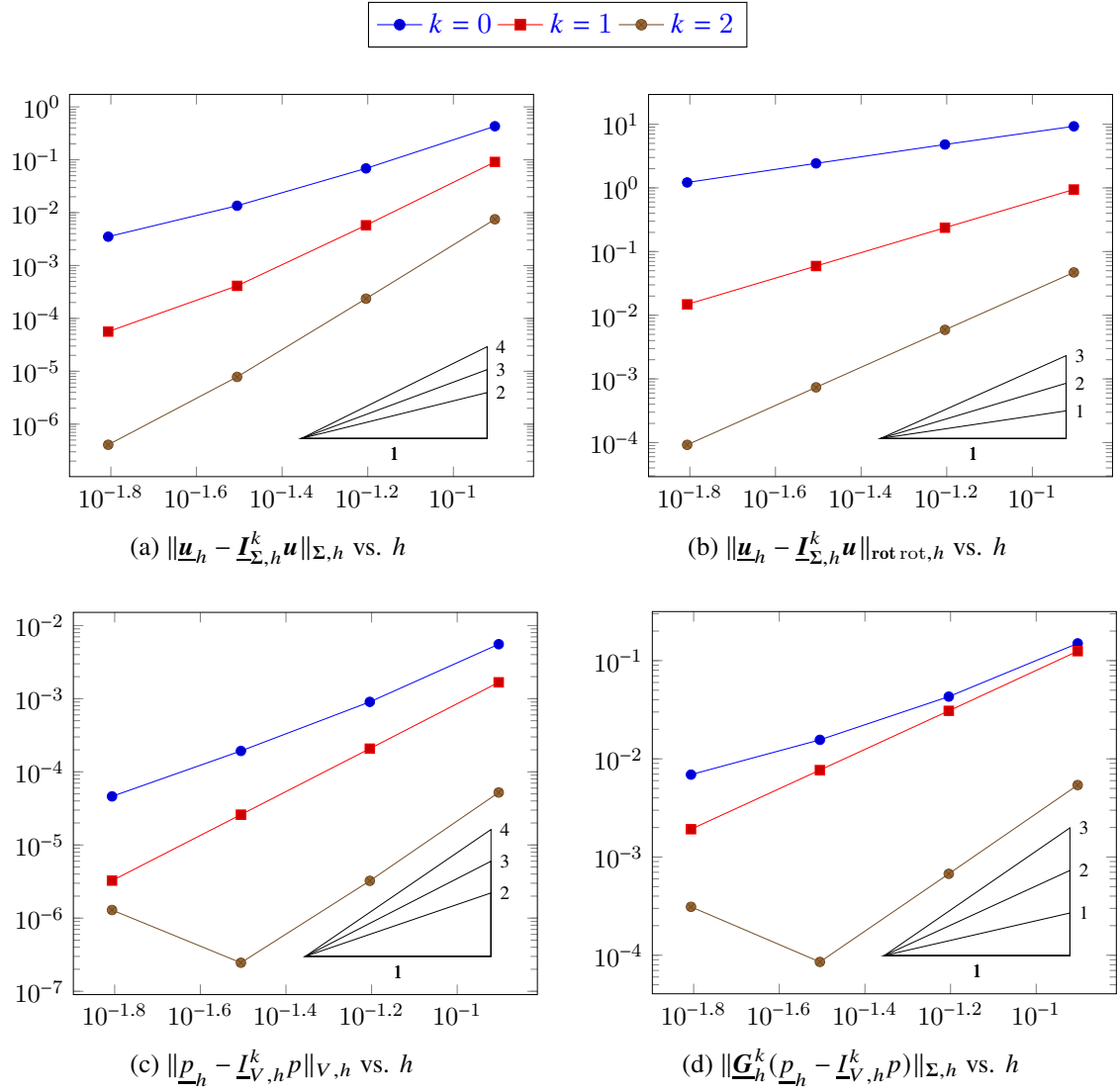

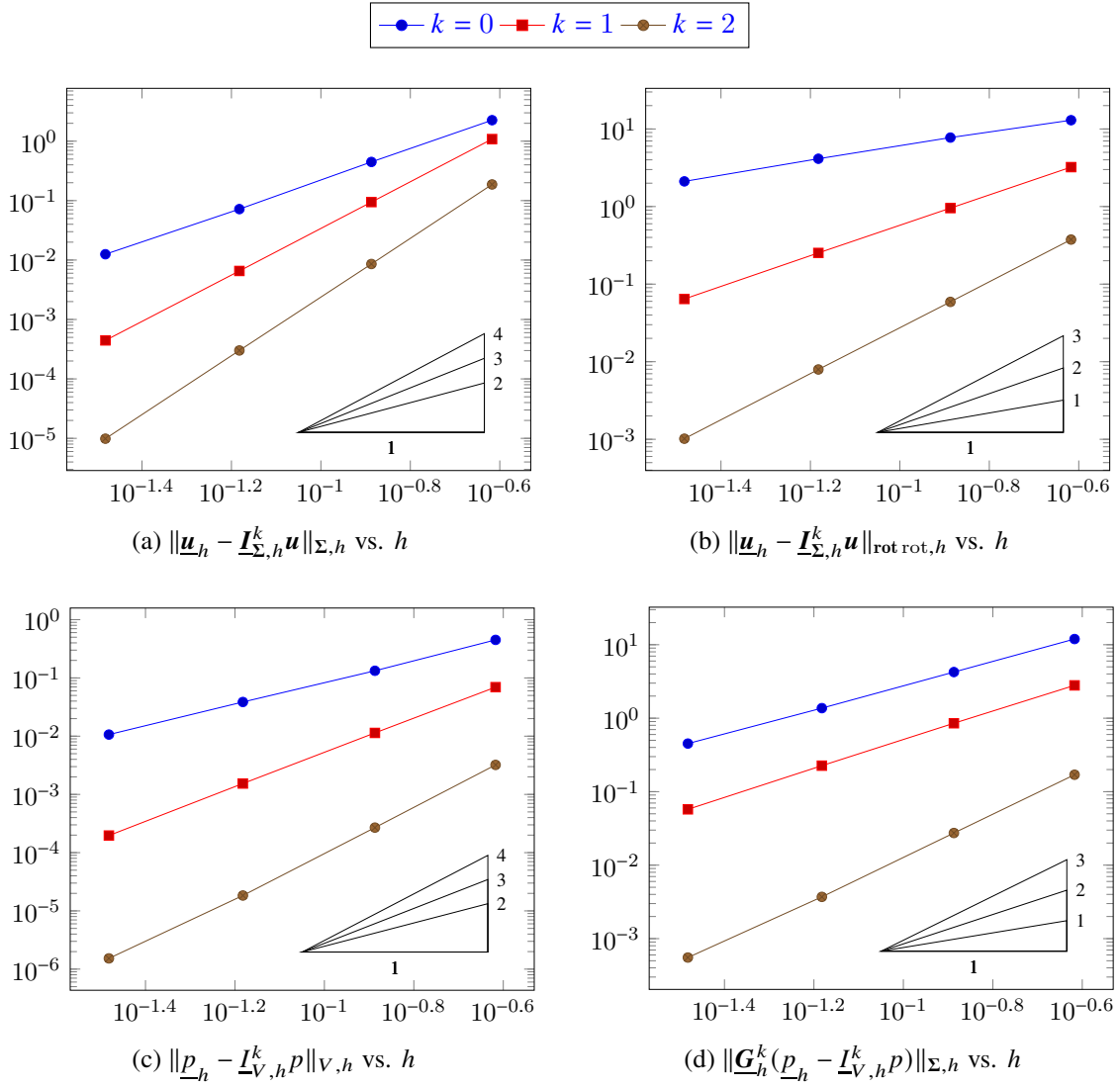
\begin{figure}\centering
  \ref{leg:hexagonal}
  \vspace{0.5cm}\\
  \begin{minipage}[b]{0.475\linewidth}\centering
    \begin{tikzpicture}[scale=0.90]
      \begin{loglogaxis}[
          legend columns=-1,
          legend to name=leg:hexagonal
        ]
        \addplot table[x=MeshSize,y=ErrUL2] {dat/hexa_k0/data_rates.dat};
        \addplot table[x=MeshSize,y=ErrUL2] {dat/hexa_k1/data_rates.dat};
        \addplot table[x=MeshSize,y=ErrUL2] {dat/hexa_k2/data_rates.dat};

        \logLogSlopeTriangle{0.9}{0.4}{0.1}{2}{black};
        \logLogSlopeTriangle{0.9}{0.4}{0.1}{3}{black};
        \logLogSlopeTriangle{0.9}{0.4}{0.1}{4}{black};
        
        \legend{$k=0$,$k=1$,$k=2$};
      \end{loglogaxis}
    \end{tikzpicture}
    \subcaption{$\norm[\bvec{\Sigma},h]{\uvec{u}_h - \ISigma{h}\bvec{u}}$ vs. $h$}
  \end{minipage}
  \begin{minipage}[b]{0.475\linewidth}\centering
    \begin{tikzpicture}[scale=0.90]
      \begin{loglogaxis}
        \addplot table[x=MeshSize,y=ErrURotRot] {dat/hexa_k0/data_rates.dat};
        \addplot table[x=MeshSize,y=ErrURotRot] {dat/hexa_k1/data_rates.dat};
        \addplot table[x=MeshSize,y=ErrURotRot] {dat/hexa_k2/data_rates.dat};

        \logLogSlopeTriangle{0.9}{0.4}{0.1}{1}{black};
        \logLogSlopeTriangle{0.9}{0.4}{0.1}{2}{black};
        \logLogSlopeTriangle{0.9}{0.4}{0.1}{3}{black};
      \end{loglogaxis}
    \end{tikzpicture}
    \subcaption{$\norm[\VROT\ROT,h]{\uvec{u}_h - \ISigma{h}\bvec{u}}$ vs. $h$}
  \end{minipage}
  \vspace{0.5cm}\\
  \begin{minipage}[b]{0.475\linewidth}\centering
    \begin{tikzpicture}[scale=0.90]
      \begin{loglogaxis}
        \addplot table[x=MeshSize,y=ErrPL2] {dat/hexa_k0/data_rates.dat};
        \addplot table[x=MeshSize,y=ErrPL2] {dat/hexa_k1/data_rates.dat};
        \addplot table[x=MeshSize,y=ErrPL2] {dat/hexa_k2/data_rates.dat};

        \logLogSlopeTriangle{0.9}{0.4}{0.1}{2}{black};
        \logLogSlopeTriangle{0.9}{0.4}{0.1}{3}{black};
        \logLogSlopeTriangle{0.9}{0.4}{0.1}{4}{black};
      \end{loglogaxis}
    \end{tikzpicture}
    \subcaption{$\norm[V,h]{\underline{p}_h - \IV{h} p}$ vs. $h$}
  \end{minipage}
  \begin{minipage}[b]{0.475\linewidth}\centering
    \begin{tikzpicture}[scale=0.90]
      \begin{loglogaxis}
        \addplot table[x=MeshSize,y=ErrPGrad] {dat/hexa_k0/data_rates.dat};
        \addplot table[x=MeshSize,y=ErrPGrad] {dat/hexa_k1/data_rates.dat};
        \addplot table[x=MeshSize,y=ErrPGrad] {dat/hexa_k2/data_rates.dat};

        \logLogSlopeTriangle{0.9}{0.4}{0.1}{1}{black};
        \logLogSlopeTriangle{0.9}{0.4}{0.1}{2}{black};
        \logLogSlopeTriangle{0.9}{0.4}{0.1}{3}{black};
      \end{loglogaxis}
    \end{tikzpicture}
        \subcaption{$\norm[\bvec{\Sigma},h]{\uGh(\underline{p}_h - \IV{h} p)}$ vs. $h$}
  \end{minipage}
  \caption{Convergence of the error components for the numerical test of Section \ref{sec:application:numerical.results} on the hexagonal mesh family.\label{fig:hexagonal}}
\end{figure}

To validate numerically the error estimates of the previous section, we solve on the unit square domain $\Omega = (0,1)^2$ the problem corresponding to the analytical solution such that, for all $\bvec{x} = \begin{pmatrix}x_1\\x_2
\end{pmatrix}\in\overline{\Omega}$,
\begin{gather*}
  \bvec{u}(\bvec{x}) = \begin{pmatrix}
    -\sin(\pi (x_1 + x_2)) \\ \sin(\pi (x_1 + x_2))
  \end{pmatrix},\qquad
  \ROT\bvec{u}(\bvec{x}) = 2\pi\cos(\pi (x_1 + x_2)),\qquad
  p(\bvec{x}) = \sin(\pi x_1) \sin(\pi x_2),
  \\
  \bvec{f}(\bvec{x}) = \begin{pmatrix}
    -4\pi^4\sin(\pi(x_1 + x_2)) + \pi\cos(\pi x_1)\sin(\pi x_2) \\
    4\pi^4\sin(\pi(x_1 + x_2)) + \pi\sin(\pi x_1)\cos(\pi x_2)
  \end{pmatrix}.
\end{gather*}
We display in Figures \ref{fig:cartesian}--\ref{fig:hexagonal} the convergence of each term in the left-hand side of the estimate \eqref{eq:err.estimate:operator.norms} with respect to the meshsize $h$ on refined sequences of Cartesian orthogonal, matching triangular, and (predominantly) hexagonal meshes.
It can be observed that the error estimate is sharp for the second and fourth component, while convergence in $h^{k+2}$ is observed for the $L^2$-like norms of the errors on $\bvec{u}$ and $p$.
A full theoretical justification of this improved $L^2$-convergence can be obtained using standard arguments based on the Aubin--Nitsche trick; for the sake of brevity, the details are omitted.
It is also worth noticing that saturation of the pressure errors is observed in Figure \ref{fig:triangular} for the finest triangular mesh and $k=2$.
This can be justified observing that the last refinement of the triangular mesh sequence corresponds to the finest mesh among the ones considered in the tests.

%------------------------------------------------------------------------------%
% Cartesian mesh family
%------------------------------------------------------------------------------%

%------------------------------------------------------------------------------%
% Triangular mesh family
%------------------------------------------------------------------------------%

%------------------------------------------------------------------------------%
% Hexagonal mesh family
%------------------------------------------------------------------------------%

\subsection{Stability analysis}\label{sec:application:stability}

\begin{proof}[Proof of Lemma \ref{lem:inf-sup}]
  Denote by $\$$ the supremum in the right-hand side of \eqref{eq:inf-sup}.
  Taking $(\uvec{v}_h, \underline{q}_h) = (\uvec{w}_h + \uGh\underline{r}_h, \underline{r}_h)$ in \eqref{eq:Ah}, noticing that $a_h(\uvec{w}_h,\uvec{w}_h) = \norm[\VROT\ROT,h]{\uvec{w}_h}^2$ (compare \eqref{eq:aT} with \eqref{eq:tnorm.rot.rot}), and recalling that $\uRh\uGh\underline{r}_h = \uvec{0}$ by the complex property proved in Theorem \ref{thm:exactness.bc}, we have
   \[
    \begin{aligned}
      \norm[\VROT\ROT,h]{\uvec{w}_h}^2
      + \cancel{a_h(\uvec{w}_h, \uGh\underline{r}_h)}
      + \norm[\bvec{\Sigma},h]{\uGh\underline{r}_h}^2
      &=
      \mathcal{A}_h((\uvec{w}_h, \underline{r}_h),(\uvec{w}_h + \uGh\underline{r}_h, \underline{r}_h))
      \\
      &\le\$\norm[\bvec{\Sigma}\times V,h]{(\uvec{w}_h + \uGh\underline{r}_h, \underline{r}_h)}
      \\
      &\lesssim\$\norm[\bvec{\Sigma}\times V,h]{(\uvec{w}_h, \underline{r}_h)},
    \end{aligned}
  \]
  where we have used the definition of the supremum to pass to the second line
  and concluded using a triangle inequality followed by the definition \eqref{eq:norm.Sigma.x.V.h} of $\norm[\bvec{\Sigma}\times V,h]{{\cdot}}$ along with $\norm[\VROT\ROT,h]{\uGh\underline{r}_h} = 0$ (again consequence of $\uRh\uGh\underline{r}_h = \uvec{0}$) to infer $\norm[\bvec{\Sigma}\times V,h]{(\uGh\underline{r}_h,\underline{0})} = \tnorm[\bvec{\Sigma},h]{\uGh\underline{r}_h}\le\tnorm[\bvec{\Sigma}\times V,h]{(\uvec{w}_h,\underline{r}_h)}$.
  Recalling the norm equivalence \eqref{eq:norm.equivalence:SSigma.h}, this yields
  \begin{equation}\label{eq:inf-sup:0}
    \norm[\VROT\ROT,h]{\uvec{w}_h}^2
    + \tnorm[\bvec{\Sigma},h]{\uGh\underline{r}_h}^2
    \lesssim\$\norm[\bvec{\Sigma}\times V,h]{(\uvec{w}_h, \underline{r}_h)}.
  \end{equation}

  The Poincar\'e inequality \eqref{eq:poincare:VZh} in $\SVZ{h}$ followed by \eqref{eq:inf-sup:0} gives
  \begin{equation}\label{eq:inf-sup:1}
    \tnorm[V,h]{\underline{r}_h}^2
    \lesssim\tnorm[\bvec{\Sigma},h]{\uGh\underline{r}_h}^2
    \lesssim\$\norm[\bvec{\Sigma}\times V,h]{(\uvec{w}_h, \underline{r}_h)}.
  \end{equation}

  It only remains to estimate $\tnorm[\bvec{\Sigma},h]{\uvec{w}_h}$.
  To this purpose, we write $\uvec{w}_h = \uvec{v}_{\bvec{w},h} + \uGh\underline{q}_{\bvec{w},h}$ with $\underline{q}_{\bvec{w},h}\in\SVZ{h}$ and $\uvec{v}_{\bvec{w},h}\in\big[\Ker\uRh\big]^\perp = \big[\Image\uGh\big]^\perp$ (recall Theorem \ref{thm:exactness.bc}), with orthogonal taken with respect to the $(\cdot,\cdot)_{\bvec{\Sigma},h}$ product.
  Notice that, since edge and vertex components of $\uGh$ not contained in $\XrotZ{h}$ vanish by definition \eqref{eq:uGh}, $\Resrot\uvec{v}_{\bvec{w},h}$ satisfies \eqref{eq:poincare:XrotZh:orthogonality} with $[\cdot,\cdot]_{\ROT,h}$ equal to the standard DDR $L^2$-product of $\Xrot{h}$ (the local version of which corresponds to the first two terms in the right-hand side of \eqref{eq:inner.prod:Sigma.T}).
  We can then write
  \begin{equation}\label{eq:inf-sup:est.w:1}
    \norm[\bvec{\Sigma},h]{\uvec{v}_{\bvec{w},h}}^2
    \stackrel{\eqref{eq:norm.equivalence:SSigma.h}}{\lesssim}\tnorm[\bvec{\Sigma},h]{\uvec{v}_{\bvec{w},h}}^2
    \stackrel{\eqref{eq:poincare:SSigmaZh:VROT}}{\lesssim}\norm[\VROT\ROT,h]{\uvec{v}_{\bvec{w},h}}^2
    = \norm[\VROT\ROT,h]{\uvec{w}_h}^2
    \stackrel{\eqref{eq:inf-sup:0}}{\lesssim}\$\norm[\bvec{\Sigma}\times V,h]{(\uvec{w}_h, \underline{r}_h)},
  \end{equation}
  where the equality follows observing that $\uRh\uvec{v}_{\bvec{w},h} = \uRh\uvec{w}_h$ (since $\uRh\uGh\underline{q}_{\bvec{w},h} = \underline{0}$) and recalling that, by definition, $\norm[\VROT\ROT,h]{\uvec{v}_{\bvec{w},h}}$ only depends on $\uvec{v}_{\bvec{w},h}$ through this quantity (see \eqref{eq:tnorm.rot.rot}).
  Taking now $(\uvec{v}_h, \underline{q}_h) = (\uvec{0}, -\underline{q}_{\bvec{w},h})$ in the expression \eqref{eq:Ah} of $\mathcal{A}_h$, we get
  \[
  \mathcal{A}_h((\uvec{w}_h, \underline{r}_h), (\uvec{0}, -\underline{q}_{\bvec{w},h}))
  = b_h(\uvec{w}_h, \underline{q}_{\bvec{w},h})
  = \cancel{(\uvec{v}_{\bvec{w},h}, \uGh\underline{q}_{\bvec{w},h})_{\bvec{\Sigma},h}}
  + \norm[\bvec{\Sigma},h]{\uGh\underline{q}_{\bvec{w},h}}^2,
  \]
  where the last equality is a consequence of the definition \eqref{eq:ah.bh.lh} of $b_h$,
  while the cancellation follows from the $(\cdot,\cdot)_{\bvec{\Sigma},h}$-orthogonality of the decomposition of $\uvec{w}_h$.
  Hence,
  \begin{equation}\label{eq:inf-sup:est.w:2}
    \norm[\bvec{\Sigma},h]{\uGh\underline{q}_{\bvec{w},h}}^2
    \lesssim\$\norm[\bvec{\Sigma}\times V,h]{(\uvec{0},\underline{q}_{\bvec{w},h})}
    \lesssim\$\norm[\bvec{\Sigma}\times V,h]{(\uvec{w},\underline{r}_h)},
  \end{equation}
  where the last equality follows observing that
  \begin{multline*}
  \norm[\bvec{\Sigma}\times V,h]{(\uvec{0},\underline{q}_{\bvec{w},h})}
  \stackrel{\eqref{eq:norm.Sigma.x.V.h}}{\le}\tnorm[V,h]{\underline{q}_{\bvec{w},h}}
  + \tnorm[\bvec{\Sigma},h]{\uGh\underline{q}_{\bvec{w},h}}  
  \stackrel{\eqref{eq:poincare:VZh}}{\lesssim}\tnorm[\bvec{\Sigma},h]{\uGh\underline{q}_{\bvec{w},h}}
  \\
  \stackrel{\eqref{eq:norm.equivalence:SSigma.h}}{\lesssim}\norm[\bvec{\Sigma},h]{\uGh\underline{q}_{\bvec{w},h}}
  \le\norm[\bvec{\Sigma},h]{\uvec{w}_h}
  \stackrel{\eqref{eq:norm.equivalence:SSigma.h}}{\lesssim}\tnorm[\bvec{\Sigma},h]{\uvec{w}_h}
  \stackrel{\eqref{eq:norm.Sigma.x.V.h}}{\le}\norm[\bvec{\Sigma}\times V,h]{(\uvec{w},\underline{r}_h)},
  \end{multline*}
  the inequality $\norm[\bvec{\Sigma},h]{\uGh\underline{q}_{\bvec{w},h}}\le\norm[\bvec{\Sigma},h]{\uvec{w}_h}$ being a consequence of the $(\cdot,\cdot)_{\bvec{\Sigma},h}$-orthogonality of the decomposition $\uvec{w}_h = \uvec{v}_{\bvec{w},h} + \uGh\underline{q}_{\bvec{w},h}$.
  We can then write  
  \begin{equation}\label{eq:inf-sup:2}
    \tnorm[\bvec{\Sigma},h]{\uvec{w}_h}^2
    \stackrel{\eqref{eq:norm.equivalence:SSigma.h}}{\lesssim}
    \norm[\bvec{\Sigma},h]{\uvec{w}_h}^2
    = \norm[\bvec{\Sigma},h]{\uvec{v}_{\bvec{w},h}}^2
    + \norm[\bvec{\Sigma},h]{\uGh\underline{q}_{\bvec{w},h}}^2
    \stackrel{\eqref{eq:inf-sup:est.w:1}, \eqref{eq:inf-sup:est.w:2}}{\lesssim}\$\norm[\bvec{\Sigma}\times V,h]{(\uvec{w},\underline{r}_h)}.
  \end{equation}
  Summing \eqref{eq:inf-sup:0}, \eqref{eq:inf-sup:1}, \eqref{eq:inf-sup:2} and simplifying, \eqref{eq:inf-sup} follows.
\end{proof}

\subsection{Convergence analysis}\label{sec:application:err.estimate}

\subsubsection{Dual consistency for rot-rot}

Given $\bvec{w}\in\bvec{\Sigma}$ (with $\bvec{\Sigma}$ denoting the domain of $\ISigma{h}$) such that $(\VROT\ROT)^2\bvec{w}\in\bvec{L}^2(\Omega;\Real^2)$, we define the rot-rot adjoint consistency error linear form $\dErotrot(\bvec{w};\cdot):\SSigmaZ{h}\to\Real$ such that, for all $\uvec{v}_h\in\SSigmaZ{h}$,
\[
\dErotrot(\bvec{w};\uvec{v}_h)
\coloneq
\ell_h((\VROT\ROT)^2\bvec{w},\uvec{v}_h)
- a_h(\ISigma{h}\bvec{w},\uvec{v}_h).
\]

\begin{lemma}[Adjoint consistency for rot-rot]\label{lem:adjoint.consistency:rot-rot}
  Recall the second complex of \eqref{eq:commutation}
  and let $\bvec{w}\in\bvec{\Sigma}$ be such that
  $\ROT\bvec{w}\in H^{k+2}(\Th)$,
  $\VROT\ROT\bvec{w}\in\bvec{H}^1(\Omega;\Real^2)\cap\bvec{H}^{k+1}(\Th;\Real^2)$,
  $\ROT(\VROT\ROT\bvec{w})\in H^1(\Omega)\cap H^{k+2}(\Th)$,
  and $(\VROT\ROT)^2\bvec{w}\in\bvec{L}^2(\Omega;\Real^2)$.
  Then, it holds
  \begin{equation}\label{eq:est:dErotrot}
    |\dErotrot(\bvec{w};\uvec{v}_h)|
    \lesssim h^{k+1}\mathcal{N}(\bvec{w})\left(
    \tnorm[\bvec{\Sigma},h]{\uvec{v}_h}
    + \norm[\VROT\ROT,h]{\uvec{v}_h}
    \right)
    \qquad\forall\uvec{v}_h\in\SSigmaZ{h},
  \end{equation}
  with $\mathcal{N}(\bvec{w})$ defined by \eqref{eq:mathcal.N}.
\end{lemma}

\begin{proof}
  Expanding first $(\cdot,\cdot)_{\bvec{\Sigma},h}$, $a_h$, and $\ell_h$ according to the respective definitions \eqref{eq:inner.prod:Sigma.h} and \eqref{eq:ah.bh.lh}, then inserting $\pm\int_T\ROT(\VROT\ROT\bvec{w})~\RT\uvec{v}_T$ inside the sum over $T\in\Th$, we obtain the following decomposition:
  \begin{equation}\label{eq:dErotrot:decomposition}
    \begin{aligned}
      \dErotrot(\bvec{w};\uvec{v}_h)
      &=
      \sum_{T\in\Th}\left[
        \int_T(\VROT\ROT)^2\bvec{w}\cdot\PSigmaT\uvec{v}_T
        - \int_T\ROT(\VROT\ROT\bvec{w})~\RT\uvec{v}_T
        \right]
      \\
      &\quad
      + \sum_{T\in\Th}\left[
        \int_T\ROT(\VROT\ROT\bvec{w})~\RT\uvec{v}_T
        - \int_T\VRT\uRT\ISigma{T}\bvec{w}\cdot\VRT\uRT\uvec{v}_T
        \right]
      \\
      &\quad    
      + \sum_{T\in\Th} s_T(\uRT\ISigma{T}\bvec{w}, \uRT\bvec{v}_T).
    \end{aligned}
  \end{equation}
  Denote by $\term_1$, $\term_2$, and $\term_3$ the terms in the first, second, and third line, respectively.
  \medskip\\
  \underline{(i) \emph{Estimate of $\term_1$.}}
  We set, for the sake of brevity, $\varphi\coloneq\ROT(\VROT\ROT\bvec{w})$ and notice that, for all $q_T\in\Poly{k+1}(T)$, it holds by \eqref{eq:PSigmaT} with $(q,\bvec{w}) = (q_T,\bvec{0})$,
  \[
  \int_T\RT\uvec{v}_T~q_T
  + \sum_{E\in\ET}\omega_{TE}\int_E v_E~q_T
  - \int_T\PSigmaT\uvec{v}_T\cdot\VROT q_T
  = 0.
  \]
  Taking $q_T = \lproj{k+1}{T}\varphi$ and inserting the above quantity inside the sum over mesh elements, we get
  \[
  \begin{aligned}
    \term_1
    &=\sum_{T\in\Th}\bigg[
    \cancel{\int_T(q_T - \varphi)~\RT\uvec{v}_T}
    + \int_T\VROT(\varphi - q_T)\cdot\PSigmaT\uvec{v}_T
    + \sum_{E\in\ET}\omega_{TE}\int_E v_E(\varphi - q_T)
    \bigg],
  \end{aligned}
  \]
  where the cancellation is a consequence of the definition of the $L^2$-orthogonal projector along with $\RT\uvec{v}_T\in\Poly{k}(T)\subset\Poly{k+1}(T)$ for all $T\in\Th$ and we have used the continuity of $\varphi$ across internal edges along with $v_E = 0$ for all $E\in\Ehb$ to insert $\varphi_{|E}$ into the boundary term.
  Using Cauchy--Schwarz inequalities on the integrals and sums, we get
  \[
  \begin{aligned}
    |\term_1|    
    &\le\left[
      \sum_{T\in\Th}\left(
      \norm[\bvec{L}^2(T;\Real^2)]{\VROT(\varphi - q_T)}^2
      + h_T^{-1} \norm[L^2(\partial T)]{\varphi - q_T}^2
      \right)
      \right]^{\nicefrac12}
    \\
    &\quad
    \times\left[
      \sum_{T\in\Th}\left(
      \norm[\bvec{L}^2(T;\Real^2)]{\PSigmaT\uvec{v}_T}^2
      + h_T\sum_{E\in\ET}\norm[L^2(E)]{v_E}^2
      \right)
      \right]^{\nicefrac12}.
  \end{aligned}
  \]
  Using the approximation properties of $\lproj{k+1}{T}$, it is readily inferred that the first factor is $\lesssim h^{k+1}\seminorm[H^{k+2}(\Th)]{\varphi}$.
  On the other hand, by the continuity property $\norm[\bvec{L}^2(T;\Real^2)]{\PSigmaT\uvec{v}_T}\lesssim\tnorm[\bvec{\Sigma},T]{\uvec{v}_T}$ (which follows from the first relation in \cite[Eq.~(4.23)]{Di-Pietro.Droniou:21*1})
  and the definition of $\tnorm[\bvec{\Sigma},h]{{\cdot}}$ (see \eqref{eq:tnorm.Sigma.T} and \eqref{eq:tnorm.rot.T}),
  the second factor is $\lesssim\tnorm[\bvec{\Sigma},h]{\uvec{v}_h}$.
  Hence, we have proved that
  \begin{equation}\label{eq:adjoint.consistency:rot-rot:T1}
    |\term_1|
    \lesssim
    h^{k+1}\seminorm[H^{k+2}(\Th)]{\ROT(\VROT\ROT\bvec{w})}
    \tnorm[\bvec{\Sigma},h]{\uvec{v}_h}.
  \end{equation}
  
  \noindent\underline{(ii) \emph{Estimate of $\term_2$.}}
  Set, for the sake of brevity $\bvec{\zeta}\coloneq\VROT\ROT\bvec{w}$.
  Letting, for all $T\in\Th$, $\bvec{z}_T\in\vPoly{k}(T;\Real^2)$ and writing the definition \eqref{eq:VRT} of $\VRT$ with $(\underline{r}_T, \bvec{v}) = (\uRT\uvec{v}_T, \bvec{z}_T)$, we infer
  \[  
  \int_T\VRT\uRT\uvec{v}_T\cdot\bvec{z}_T
  - \int_T\RT\uvec{v}_T~\ROT\bvec{z}_T
  - \sum_{E\in\ET}\omega_{TE}\int_E\PWE\uRE\uvec{v}_E~(\bvec{z}_T\cdot\tangent_E)
  = 0.
  \]
  Adding this quantity to $\term_2$
  and using the continuity of the tangential trance of $\bvec{\zeta}$ at interfaces along with the fact that $\PWE\uRE\uvec{v}_E = 0$ for all $E\in\Ehb$ (since $\uvec{v}_h\in\SSigmaZ{h}$) to insert $\bvec{\zeta}_{|E}\cdot\tangent_E$ it into the boundary integral,
  we get
  \begin{multline*}
    \term_2 = \sum_{T\in\Th}\bigg[
      \int_T\left(\bvec{z}_T - \VRT\IW{T}\ROT\bvec{w}\right)\cdot\VRT\uRT\uvec{v}_T
      \\
      + \int_T\ROT(\bvec{\zeta} - \bvec{z}_T)~\RT\uvec{v}_T
      + \sum_{E\in\ET}\omega_{TE}\int_E(\bvec{\zeta} - \bvec{z}_T)\cdot\tangent_E~\PWE\uRE\uvec{v}_E
      \bigg],
  \end{multline*}
  where we have additionally used the commutation property expressed by the rightmost portion of \eqref{eq:commutation} to replace $\RT\ISigma{T}\bvec{w}$ with $\IW{T}\ROT\bvec{w}$ in the first term.
  Take $\bvec{z}_T = \vlproj{k}{T}\bvec{\zeta}$ for all $T\in\Th$.
  Integrating by parts the second term in square brackets, we obtain
  \begin{multline*}
    \term_2 = \sum_{T\in\Th}\bigg[
      \int_T\left(\bvec{z}_T - \VRT\IW{T}\ROT\bvec{w}\right)\cdot\VRT\uRT\uvec{v}_T
      \\
      + \cancel{\int_T(\bvec{\zeta} - \bvec{z}_T)\cdot\VROT\RT\uvec{v}_T}
      + \sum_{E\in\ET}\omega_{TE}\int_E(\bvec{\zeta} - \bvec{z}_T)\cdot\tangent_E~(\PWE\uRE\uvec{v}_E - \RT\uvec{v}_T)
      \bigg],
  \end{multline*}
  where the cancellation follows by definition of $\vlproj{k}{T}$ after observing that $\VROT\RT\uvec{v}_T\in\vPoly{k-1}(T;\Real^2)\subset\vPoly{k}(T;\Real^2)$.
  Using Cauchy--Schwarz inequalities on the integrals and the sums, we obtain
  \[
  \begin{aligned}
    |\term_2|
    &\lesssim\left[
      \sum_{T\in\Th}\left(
      \norm[\bvec{L}^2(T;\Real^2)]{\bvec{z}_T - \VRT\IW{T}\ROT\bvec{w}}^2
      + h_T\norm[\bvec{L}^2(\partial T;\Real^2)]{\bvec{\zeta} - \bvec{z}_T}^2
      \right)
      \right]^{\nicefrac12}
    \\
    &\quad\times
    \left[
    \sum_{T\in\Th}\left(
    \norm[\bvec{L}^2(T;\Real^2)]{\VRT\uRT\uvec{v}_T}^2
    + h_T^{-1}\hspace{-1ex}\sum_{E\in\ET}\norm[L^2(E)]{\RT\uvec{v}_T - \PWE\uRE\uvec{v}_E}^2
    \right)
  \right]^{\nicefrac12}\hspace{-1ex}.
  \end{aligned}
  \]
   Using the approximation properties of $\VRT\IW{T}$ (which can be proved using similar arguments as for \cite[Eq. (6.5)]{Di-Pietro.Droniou:21*1}) and $\vlproj{k}{T}$, it can be checked that the first factor in the right-hand side is $\lesssim h^{k+1}\seminorm[\bvec{H}^{k+1}(\Th;\Real^2)]{\bvec{\zeta}}$.
  By definition of $\norm[\VROT\ROT,h]{{\cdot}}$ along with \eqref{eq:poincare:rot-rot:estimate.rhs}, the second factor is $\lesssim\norm[\VROT\ROT,h]{\uvec{v}_h}$.
  Plugging the above estimates into the bound of $|\term_2|$, we arrive at
  \begin{equation}\label{eq:adjoint.consistency:rot-rot:T2}
    |\term_2|\lesssim h^{k+1}\seminorm[\bvec{H}^{k+1}(\Th;\Real^2)]{\VROT\ROT\bvec{w}}\norm[\VROT\ROT,h]{\uvec{v}_h}.
  \end{equation}
  \noindent\underline{(iii) \emph{Estimate of $\term_3$.}} Moving to the third and last term, we use again the fact that $\uRT\ISigma{T}\bvec{w} = \IW{T}\ROT\bvec{w}$ by the rightmost commutative diagram of \eqref{eq:commutation}, then use Cauchy--Schwarz inequalities  to write
  \[
  |\term_3|
  \lesssim h^{k+1} \left(
  \sum_{T\in\Th} s_T(\IW{T}\ROT\bvec{w},\IW{T}\ROT\bvec{w})
  \right)^{\nicefrac12}
  \left(
  \sum_{T\in\Th} s_T(\uRT\uvec{v}_T, \uRT\uvec{v}_T)
  \right)^{\nicefrac12}.
  \]
  Comparing \eqref{eq:tnorm.rot.rot} and \eqref{eq:sT}, it is clear that the second factor is $\lesssim\norm[\VROT\ROT,h]{\uvec{v}_h}$.
  In order to estimate the first factor, we start by recalling the following boundedness property for $\IW{T}$, which is the two-dimensional counterpart of \cite[Eq. (4.26)]{Di-Pietro.Droniou:21*1} and can be proved using trace inequalities:
  For all $r\in H^2(T)$,
  \begin{equation}\label{eq:IWT:boundedness}
    \tnorm[W,T]{\IW{T}r}
    \lesssim\norm[L^2(T)]{r}
    + h_T\seminorm[H^1(T)]{r}
    + h_T^2\seminorm[H^2(T)]{r}.
  \end{equation}
  Assuming that we can prove, for all $T\in\Th$ and all $\underline{r}_T\in\SW{T}$, that
  \begin{equation}\label{eq:boundedness:sT}
    s_T(\underline{r}_T,\underline{r}_T)^{\nicefrac12}
    \lesssim h_T^{-2}\tnorm[W,T]{\IW{T}\PWT\underline{r}_T - \underline{r}_T}^2,
  \end{equation}
  we can then write
  \[
  \begin{aligned}
    s_T(\IW{T}\ROT\bvec{w},\IW{T}\ROT\bvec{w})^{\nicefrac12}
    &\!\stackrel{\eqref{eq:boundedness:sT}}{\lesssim}
    h_T^{-2}\tnorm[W,T]{\IW{T}(\PWT\IW{T}\ROT\bvec{w} - \ROT\bvec{w})}
    \\
    &
    \!\stackrel{\eqref{eq:IWT:boundedness}}{\lesssim}
    h_T^{-2}\norm[L^2(T)]{\PWT\IW{T}\ROT\bvec{w} - \ROT\bvec{w}}
    \\
    &\quad
    + h_T^{-1}\seminorm[H^1(T)]{\PWT\IW{T}\ROT\bvec{w} - \ROT\bvec{w}}
    \\
    &\quad
    + h_T\seminorm[H^2(T)]{\PWT\IW{T}\ROT\bvec{w} - \ROT\bvec{w}}
    \\
    &\lesssim
    h_T^{k+1}\seminorm[H^{k+2}(T)]{\ROT\bvec{w}},
  \end{aligned}
  \]
  where we have used the approximation properties of $\PWT$ (which are the two-dimensional counterpart of \cite[Eq. (6.2)]{Di-Pietro.Droniou:21*1}) to conclude.
  This estimate leads to 
  \begin{equation}\label{eq:adjoint.consistency:rot-rot:T3}
    |\term_3|\lesssim h^{k+1}\seminorm[H^{k+2}(\Th)]{\ROT\bvec{w}}\norm[\VROT\ROT,h]{\uvec{v}_T}.
  \end{equation}
  It only remains to prove \eqref{eq:boundedness:sT}.
  Expanding $s_T(\underline{r}_T,\underline{r}_T)$ according to its definition \eqref{eq:sT}, we have
  \[
  \begin{aligned}
    s_T(\underline{r}_T, \underline{r}_T)
    &= h_T^{-2}\norm[L^2(T)]{\lproj{k}{T}(\PWT\underline{r}_T - r_T)}^2
    + h_T^{-1}\sum_{E\in\ET}\norm[L^2(E)]{\PWT\underline{r}_T - \PWE\underline{r}_E}^2
    \\
    &\!\stackrel{\eqref{eq:norm.equivalence:1d}}{\lesssim}
    h_T^{-2}\norm[L^2(T)]{\lproj{k}{T}(\PWT\underline{r}_T - r_T)}^2
    \\
    &\quad
    + h_T^{-1}\sum_{E\in\ET}
    \norm[L^2(E)]{\lproj{k-1}{E}\PWT\underline{r}_T - \underbrace{\lproj{k-1}{E}\PWE\underline{r}_E}_{\stackrel{\eqref{eq:PWE}}{=}r_E}}^2
    \\
    &\quad
    + h_T^{-1}\sum_{E\in\ET}h_E\sum_{\nu\in\VE}\big(
    \PWT\underline{r}_T(\bvec{x}_\nu) - \underbrace{\PWE\underline{r}_E(\bvec{x}_\nu)}_{\stackrel{\eqref{eq:PWE}}{=}r_\nu}
    \big)^2
    \\
    &\!\stackrel{\eqref{eq:tnorm.W.h}}{\lesssim}
    h_T^{-2}\tnorm[W,T]{\IW{T}\PWT\underline{r}_T - \underline{r}_T}^2,
  \end{aligned}
  \]
  where the conclusion follows using mesh regularity to write $h_T^{-1}h_E\lesssim 1$ and noticing that $\sum_{E\in\ET}\sum_{\nu\in\VE}\alpha_\nu = 2\sum_{\nu\in\VT}\alpha_\nu$ for all families $(\alpha_\nu)_{\nu\in\VT}\in\Real^{\VT}$ for the last term.
  This concludes the proof of \eqref{eq:boundedness:sT}.
  \medskip\\
  \underline{(iv) \emph{Conclusion}.}
  Plugging \eqref{eq:adjoint.consistency:rot-rot:T1}, \eqref{eq:adjoint.consistency:rot-rot:T2}, and \eqref{eq:adjoint.consistency:rot-rot:T3} into \eqref{eq:dErotrot:decomposition} yields the conclusion.
\end{proof}

\subsubsection{Proof of Theorem \ref{thm:error.estimate}}

\begin{proof}[Proof of Theorem \ref{thm:error.estimate}]
  We use the abstract framework of \cite{Di-Pietro.Droniou:18}.
  The consistency error is
  \begin{equation}\label{eq:Err:decomposition}
    \begin{aligned}
      \Err((\bvec{u},p); (\uvec{v}_h, \underline{q}_h))
      &\coloneq\ell_h(\bvec{f}, \uvec{v}_h)
      - \mathcal{A}_h((\ISigma{h}\bvec{u}, \IV{h}p), (\uvec{v}_h, \underline{q}_h))
      \\
      &= \dErotrot(\bvec{u}; \uvec{v}_h)
      + \Egrad(p; \uvec{v}_h)
      + \dEgrad(\bvec{u}; \underline{q}_h),
    \end{aligned}
  \end{equation}
  where we have used the fact that $\bvec{f} = (\VROT\ROT)^2\bvec{u} + \GRAD p$ almost everywhere in $\Omega$ and introduced the primal and adjoint gradient error linear forms $\Egrad(p;\cdot):\SSigmaZ{h}\to\Real$ and $\dEgrad(\bvec{u}; \cdot):\SVZ{h}\to\Real$ such that, for all $(\uvec{v}_h, \underline{q}_h)\in\SSigmaZ{h}\times\SVZ{h}$,
  \[
  \begin{aligned}
    \Egrad(p; \uvec{v}_h)
    &\coloneq
    \sum_{T\in\Th}\int_T\GRAD p\cdot\PSigmaT\uvec{v}_T
    - b_h(\uvec{v}_h, \IV{h}p),
    \\
    \dEgrad(\bvec{u}; \underline{q}_h)
    &\coloneq
    \sum_{T\in\Th}\int_T\underbrace{\DIV\bvec{u}}_{=0}~\PVT\underline{q}_T
    + b_h(\ISigma{h}\bvec{u}, \underline{q}_h).
  \end{aligned}
  \]
  Using arguments analogous to the ones of \cite[Corollary~2 and Theorem~9]{Di-Pietro.Droniou:21*1}, it can be proved that
  \begin{equation}\label{eq:est:Egrad.dEgrad}
    \begin{aligned}
      |\Egrad(p; \uvec{v}_h)|
      &\lesssim h^{k+1}\seminorm[H^{k+2}(\Th)]{p}\tnorm[\bvec{\Sigma},h]{\uvec{v}_h},
      \\
      |\dEgrad(\bvec{u}; \underline{q}_h)|
      &\lesssim h^{k+1}\seminorm[\bvec{H}^{(k+1,2)}(\Th;\Real^2)]{\bvec{u}}
      \tnorm[\bvec{\Sigma},h]{\uGh\underline{q}_h}.
    \end{aligned}
  \end{equation}
  Plugging \eqref{eq:est:dErotrot} with $\bvec{w} = \bvec{u}$ and \eqref{eq:est:Egrad.dEgrad} into \eqref{eq:Err:decomposition} and recalling \eqref{eq:norm.Sigma.x.V.h}, we obtain
  \[
  |\Err((\bvec{u},p); (\uvec{v}_h, \underline{q}_h))|
  \lesssim h^{k+1}\left(
  \mathcal{N}(\bvec{u})
  + \seminorm[\bvec{H}^{(k+1,2)}(\Th;\Real^2)]{\bvec{u}}
  + \seminorm[H^{k+2}(\Th)]{p}
  \right)\norm[\bvec{\Sigma},h]{(\uvec{v}_h,\underline{q}_h)}.
  \]
  Recalling \cite[Theorem 10]{Di-Pietro.Droniou:18} concludes the proof.
\end{proof}

%------------------------------------------------------------------------------%
% Bibliography
%------------------------------------------------------------------------------%

\printbibliography

\end{document}